\documentclass[a4paper,10pt]{article}
\usepackage{stmaryrd}
\usepackage{amsfonts}
\usepackage{bbm}
\usepackage{amscd}
\usepackage{mathrsfs}
\usepackage{latexsym,amssymb,amsmath,amscd,amscd,amsthm,amsxtra}
\usepackage[dvips]{graphicx}
\usepackage[utf8]{inputenc}
\usepackage[T1]{fontenc}
\usepackage{lmodern}
\usepackage{amssymb}
\usepackage[all]{xy}
\usepackage{nicefrac,mathtools,enumitem}
\usepackage{microtype}
\usepackage{CJK}
\usepackage{amssymb}
\usepackage{epstopdf}
\usepackage{graphicx}
\usepackage{graphics}
\usepackage[T1]{fontenc}
\usepackage[utf8]{inputenc}
\usepackage{authblk}
\usepackage{subfigure}
\usepackage{caption}
\usepackage{framed}
\usepackage{mathrsfs}
\usepackage[noend]{algpseudocode}
\usepackage{algorithmicx,algorithm}
\usepackage{framed}
\usepackage{booktabs}
\usepackage{threeparttable}

\usepackage{diagbox}
\usepackage{graphicx}
\usepackage{float}
\usepackage{subfigure}

\numberwithin{equation}{section}
\newtheorem{thm}{Theorem}
\newtheorem{lem}{Lemma}
\newtheorem{cor}{Corollary}

\newtheorem{rmk}{Remark}

\newtheorem{ass}{Assumption}

\newcommand {\emptycomment}[1]{}

\newcommand{\be }{\begin{equation}}
\newcommand{\ee }{\end{equation}}

\def\bea{\begin{eqnarray}}
\def\eea{\end{eqnarray}}
\def\be{\begin{equation}}
\def\ee{\end{equation}}
\def\blm{\begin{lem}}
\def\elm{\end{lem}}

\begin{document}


\title{Distributed Constrained  Optimization with Delayed Subgradient Information over Time-Varying Network under Adaptive Quantization\thanks{Jie Liu, Zhan Yu, Daniel W. C. Ho are with the Department of Mathematics, City University of Hong Kong, Hong Kong. Email: jliu285-c@my.cityu.edu.hk; zhanyu2-c@my.cityu.edu.hk; madaniel@cityu.edu.hk}}
\author{Jie Liu, Zhan Yu, Daniel W. C. Ho}
\date{}
\maketitle

\begin{abstract}

In this paper, we consider a distributed constrained optimization problem with delayed subgradient information over the time-varying communication network, where each agent can only communicate with its neighbors and the communication channel has a limited data rate. We propose an adaptive quantization method to address this problem. A mirror descent algorithm with delayed subgradient information is established based on the theory of Bregman divergence. With non-Euclidean Bregman projection-based scheme, the proposed method essentially generalizes many previous classical Euclidean projection-based distributed algorithms. Through the proposed adaptive quantization method, the optimal value without any quantization error can be obtained. Furthermore, comprehensive analysis on convergence of the algorithm is carried out  and our results show that the optimal convergence rate $O(1/\sqrt{T})$ can be obtained under appropriate conditions. Finally, numerical examples are presented to demonstrate the effectiveness of our algorithm and theoretical results.

\textbf{Keywords:} Distributed Optimization, Mirror Descent Algorithm, Adaptive Quantization, Delayed Subgradient Information, Multiagent Network.
\end{abstract}

\section{Introduction}
Recently, the distributed  optimization algorithms for the network system have been studied widely (see in \cite{Nedic1}-\cite{Yu2}). In the distributed optimization problem, there is no central coordination between different agents. Each agent knows about its local function and can only communicate with its neighboring agents in the network.
The objective function is composed of sum of local functions. Additionally, these agents, by sending updated information to their neighboring agents, cooperatively minimize the objective function. These distributed methods are critical in many engineering problems, such as localization in sensor networks   \cite{Huangchi}, smart grid optimization  \cite{Deng000}, aggregative games  \cite{Deng00},   resource allocation \cite{YI0,Deng0}, decentralized estimation \cite{Chenbo} and distributed control problems \cite{Song}.

The purpose of distributed optimization algorithms is to solve  optimization problems through distributed process in which the agents cooperatively minimize the objective function via information communication. The information communication is often carried out between an agent and its neighbours. Different kinds of distributed optimization algorithms have been proposed in recent years. In \cite{Nedic1}, the authors propose the subgradient algorithm to solve not necessarily smooth distributed optimization over time-varying communication network topology. In \cite{Nedic2}, the authors propose the distributed stochastic gradient push algorithm over time-varying directed graphs for optimization.  In \cite{Duchi1}, the authors propose distributed dual averaging algorithm  and analyze its  convergence rate. In \cite{YI1}, the authors propose distributed gradient algorithm for constrained optimization. In \cite{Qingshanliu}, the authors propose a collaborative neurodynamic approach to  distributed constrained optimization. In \cite{Shaofu1}, the authors propose a collaborative neurodynamic approach to multiple-objective distributed optimization. In  \cite{Qingshanliu2}, the authors propose a one-layer projection neural network for solving nonsmooth optimization problems with linear equalities and bound constraints. In \cite{Deng1}, the authors propose the distributed optimization for multiple heterogeneous Euler-Lagrangian systems.

Mirror descent methods for the distributed optimization have attracted much attention since it was proposed. Compared with other distributed subgradient projection methods, mirror descent algorithm uses customized Bregman divergence rather than Euclidean distance, which can be viewed as non-Euclidean projection method and a generalization of the distributed gradient method based Euclidean projection methods. Mirror descent algorithms have been studied extensively in recent year. Here are some excellent mirror descent algorithm references (see \cite{Nedic3, Li, Yuan2,Yuan5,Yu2}).  In \cite{Nedic3}, the authors propose stochastic subgradient mirror descent algorithm for constrained distributed optimization problems. In \cite{Li}, the authors propose the distributed mirror descent algorithm  over time varying multi-agent network with  delayed gradient for convex optimization. In \cite{Yuan2}, the authors propose  distributed mirror descent algorithm over  time-varying  network  consisting of multiple interacting nodes for online composite optimization. In \cite{Yuan5}, the authors propose distributed stochastic  mirror descent method over a class of time-varying network for strongly convex objective functions. In \cite{Yu2}, the authors propose distributed randomized gradient-free mirror descent algorithm for constrained optimization.

In most cases, limited communication channel between  different agents is very common and the exchanged information  need to be quantized to meet the limited communication data rate (see \cite{Pu1,Yi,Yuanquantization,Lululi}). Thus one needs to design appropriate quantizer and develops the distributed algorithms to solve distributed optimization  under limited communication channel. Previous works have proposed many distributed optimization algorithms with different kinds of quantization methods and analyze quantization's effect to the convergence (see \cite{Nedic4,Pu1,Yi}). In \cite{Nedic4},  the authors propose subgradient method under the quantization method, whose quantization value is the integer multiples of a given value, and analyze quantization error's effect to the convergence. In \cite{Pu1}, the authors propose a progressing quantization method in distributed optimization. In \cite{Yi},  the authors propose quantization method with encoder-decoder scheme and  zooming-in technique, under which the optimal value can be obtained by distributed quantized subgradient algorithm over time-varying communication network. The progressing quantization methods are  useful and significant. In \cite{Pu1}, the optimal value without quantization error can be obtained through progressing quantization.  However, the static quantizers in \cite{Nedic4,Yuanquantization} cause the non-convergence due to quantization error at each iteration. Therefore, we hope to design an adaptive quantization approach to obtain accurate optimal value. This motivates us to propose mirror descent algorithm under adaptive quantization to solve distributed optimization problem.

Time delays are unavoidable in  real life and many researchers have investigated time delay's effect to different  kinds of networks' different  status (see \cite{LU1,LU2,Yang}), such as pinning impulsive stabilization of nonlinear dynamical networks \cite{LU1}, consensus over directed static networks \cite{LU2} and synchronization of randomly coupled neural networks \cite{Yang}.

In the  distributed optimization (see \cite{Nedic1,Nedic2,Nedic3,Duchi1,Yuan5,Yu2,Yuanquantization,guanghongyang,Yuan1,Yuan9,Yuan10}), each agent needs to update parameter and calculate (sub)gradient based on local parameter in parallel.  Then each agent receives  current (sub)gradient information. However, the asynchronous process of updating parameter and calculating (sub)gradient will cause time delay. Under that circumstance, each agent receives outdated (sub)gradient information. Therefore, it is significant to study the distributed optimization algorithms with the presence of time delay. Our paper considers the asynchronous subgradient methods, where each agent  receives outdated rather than current subgradient information. By this way, each agent can update the parameters   and  compute (sub)gradients asynchronously.  The asynchrony of two processes updating parameters and calculating subgradient is very common in  real life, such as master worker architectures for distributed computation \cite{Tsianos2} and other similar model \cite{Langford,Nedic55}. Many researchers have investigated asynchronous process of distributed optimization algorithms and take delayed (sub)gradient into consideration. In \cite{Nedic55}, authors propose distributed asynchronous incremental subgradient methods with the presence of time delay. In \cite{Li}, authors propose distributed mirror descent algorithm for multi-agent optimization with delayed gradient. In \cite{Tsianos2}, authors propose  gradient-based optimization algorithms with delayed stochastic gradient information. In \cite{Wang1}, authors generalized dual averaging subgradient algorithm under delayed subgradient information.



The contribution of this paper is summerized as follows. (i) Firstly,  the distributed mirror descent algorithm with adaptive quantization is proposed to address limited communication channel. The traditional uniform quantizer uses the static quantization parameters mid-value and  interval size while  the proposed  adaptive quantizer  changes the quantization parameters mid-value and  interval size  at each iteration. The existing works such as \cite{Nedic4}, \cite{Pu1}, \cite{Yi},  \cite{Doan} design  quantizers to address the limited communication channel in distributed optimization algorithm with Euclidean projection method. However, these quantization schemes can not be directly established in the  non-Euclidean projection based methods. In this paper, we overcome this difficulty and establish the proposed distributed non-Euclidean quantization method by employing some new techniques on handling mirror descent structure. The appropriate adaptive quantizer is designed to realize the quantization in distributed optimization with non-Euclidean projection method. The proposed adaptive quantizer helps to asymptotically alleviate the quantized error but only uses a finite number of bits for quantization.

(ii)  We analyze the convergence of the mirror descent algorithm under adaptive quantization and also derive some sufficient conditions on stepsize and quantization parameter for the convergence of the proposed algorithm.  The convergence rates are comprehensively investigated by considering different stepsizes and quantization parameters. Compared with  \cite{Doan},  the convergence rate of distributed optimization algorithm in \cite{Doan} is $O(\ln T/\sqrt{T})$ while the convergence rate in this paper is $O(1/\sqrt{T})$. Our algorithm's convergence rate is faster than that in \cite{Doan}. Also, the communication network in \cite{Doan} is static while we consider a class of time-varying communication network, which is more realistic.  Compared with algorithm in \cite{Li}, an adaptive quantization method has been designed for mirror descent algorithm to address limited communication capacity and the convergence rate is still $O(1/\sqrt{T})$, which is the same with  that in \cite{Li}. Furthermore, the assumptions in this paper are much easier to be satisfied. In addition, the  objection function's subgradient has upper bound while the objection function's gradient in \cite{Li} should satisfy Lipschitz continuous.

(iii) The third contribution is that  the asynchronous operation of optimization algorithms with the presence of time delay has been considered. After careful analysis, we conclude that the convergence is guaranteed under appropriate conditions with any time delay. Finally, the algorithm can asymptotically converge to an optimal solution without quantization error.  In this paper, we significantly improve our previous works \cite{Yuan2, Yu2, Yuan5} on distributed mirror descent methods in several aspects. To the best of our knowledge, this is the first work to propose the adaptive quantization method to address limited communication channel in mirror descent algorithm and simultaneously take  delayed subgradient information into consideration in the study of distributed mirror descent method.




The rest of this paper is organized as follows. Section \uppercase\expandafter{\romannumeral2} introduces some notations and definitions. Section \uppercase\expandafter{\romannumeral3} defines the problem and propose some assumptions. Section \uppercase\expandafter{\romannumeral4} proposes the mirror descent algorithm with delayed subgradient information under adaptive quantization. Section \uppercase\expandafter{\romannumeral5} analyzes the convergence of the algorithm and discusses how to select stepsize and quantization parameter. Then we show the convergence rate of different stepsize and quantization parameter. Section \uppercase\expandafter{\romannumeral6} provides numerical simulations to verify theoretical results and Section \uppercase\expandafter{\romannumeral7} concludes this paper.

\section{Notation and Definition}

\subsection{Notation}
We first introduce some notations.
For a vector $x\in \mathbb{R}^n$, $||x||$ and $||x||_{\infty}$ are  Euclidean norm and infinity norm, respectively. The $j$th entry of  vector $x\in \mathbb{R}^n$ is $[x]_j$ and the $i$th row, $j$th column of matrix $P\in\mathbb{R}^{n\times n}$ is $[P]_{ij}$. For the $a,b,c\in \mathbb{R}^n$, we use $a\in [b,c]$ or $b\preceq a\preceq c$ to denote $[a]_j\in [[b]_j,[c]_j]$ with $j=1,\cdots,n$. The column vector $\textbf{1}\in \mathbb{R}^{n\times1}$ whose each entry is 1. For a given compact convex set $\mathcal{X}$, we use $Proj_{\mathcal{X}}\{a\}$ to denote the projection of $a$ to set $\mathcal{X}$, with
 \begin{eqnarray*}
Proj_{\mathcal{X}}\{a\}=\arg\min_{b\in\mathcal{X}}||b-a||.
\end{eqnarray*}

For a given non-smooth convex function $h(x):\mathbb{R}^n\rightarrow \mathbb{R}$, we use $\partial h(a)=\{g\in \mathbb{R}^n|h(b)\geq h(a)+g^T(b-a)\}$ is the set of the subgradient of $h(x)$ at $a$. The set $\partial h(a)$ is nonempty since $h(x)$ is convex function.

Let $h(x): A\rightarrow B$ is a $\sigma_h$ strongly convex function if and only if $h(b)\geq h(a)+\langle \nabla h(a), b-a\rangle+\frac{\sigma_{h}}{2}||b-a||^2$ for any $a,b\in A$. $L(g)$ denotes a Lipschitz constant of the function $g$ if and only if $||g(a)-g(b)||\leq L||a-b||$ for any $a,b\in$ $dom$ $g$.

For the two function $F(t)$ and $G(t)$, $F(t)=O(G(t))$ means if there are positive $T>0$ and  $C>0$ such that $F(t)\leq CG(t)$ when $t>T$.

\subsection{Uniform Quantization}

\begin{figure}[H]
\centering
\includegraphics[width=2.5in]{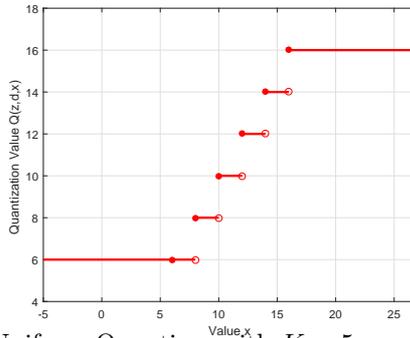}
\caption{Uniform Quantizer with $K=5$, $z=11$ and $d=5$}
\label{Quantizer}
\end{figure}

The uniform quantizer in $\mathbb{R}$ with  mid-value  $z\in \mathbb{R}$,  quantization interval size $d\in \mathbb{R}$ and a fixed number of bits $K+1$ is defined as
$$ Q(z,d,x)=\left\{
\begin{array}{rcl}
z-d,\quad\quad\quad& &{x-z\in(-\infty,-d)}\\
z+\frac{(2j-K)d}{K},& &{x-z+d\in[\frac{2 j d}{K},\frac{2(j+1)d}{K})}\\
z+d,\quad\quad\quad& &{x-z\in[d,+\infty)}
\end{array} \right. $$
where $j=0,1,\cdots,K-1$. Fig. \ref{Quantizer} is a unform quantizer with $K=5$, $z=11$ and $d=5$.

$\\$     For the quantizer in $\mathbb{R}^n$, we can also define a quantization function $\widehat{Q}(z,d,x):$ $\mathbb{R}^n\times \mathbb{R}^n\times \mathbb{R}^n\rightarrow \mathbb{R}^n$ as
follow
\begin{eqnarray}\label{quantization}
[\widehat{Q}(z,d,x)]_i=Q([z]_i,[d]_i,[x]_i).
\end{eqnarray}

\noindent where $i=1,\cdots, n$. $[\widehat{Q}(z,d,x)]_i$, $[z]_i$, $[d]_i$ and $[x]_i$ are the $i$th  element of vector $\widehat{Q}(z,d,x)$, $z$, $d$ and $x$, respectively. For the  uniform quantizer $\widehat{Q}(z,d,x)$, $z\in \mathbb{R}^n$ and $d\in \mathbb{R}^n$ denote the mid value vector  and  quantization interval size vector, respectively. When the vector $x\in \mathbb{R}^n$ falls inside the quantization interval $[z-d, z+d]$, the quantization error is bounded by
\begin{eqnarray}\label{quantizationbound}
||x-\widehat{Q}(z,d,x)||\leq \sqrt{n}||d||_{\infty}\leq \sqrt{n}||d||.
\end{eqnarray}


\section{Distributed Optimization over Time Varying Network}

In this paper, we consider an distributed optimization problem defined over time-varying communication network with $N$ nodes. The function $f_j(x)$: $\mathcal{X}\rightarrow \mathbb{R}$ with $j=1,2,\cdots,N$ are non-smooth convex functions and $\mathcal{X}\in\mathbb{R}^n$ is non-empty, convex and compact. The objective function is
\begin{eqnarray}\label{goal}
\min_{x\in\mathcal{X}}\  f(x)=\sum^N_{j=1}f_j(x).
\end{eqnarray}
The optimal value set $\mathcal{X}^*=\arg\min_{x\in\mathcal{X}} f(x)$ of problem (\ref{goal}) is not empty. There is no central coordination between the agents and each agent $j$ only knows its local function $f_j(x)$. The direct graph $\mathcal{G}(t)=\{\mathcal{V},E(t),P(t)\}$ denotes time-varying communication network topology, where $\mathcal{V}=\{1,2,\cdots, N\}$, $E(t)=\{(j,i)|$ agent $j$ and $i$ are connected, $i,j\in \mathcal{V}\}$ and $P(t)\in\mathbb{R}^{n\times n}$ is the correspond weigh matrix at the time $t$. We define $N_i(t)=\{j\in\mathcal{V}|$ $(j,i)\in E(t)\}$ as  agent $i$'s neighbor set at time $t$. The node $i$ sends information to node $j$ at time $t$ if and only if $j\in N_i(t)$. Agent $i$ and $j$ are not connected at time $t$ if and only if  $[P(t)]_{ij}=0$.

\begin{ass}\label{ass1}
The time varying network's corresponding communication matrix $P[t]$ is doubly stochastic at each time $t$, i.e. $\sum^N_{i=1}[P(t)]_{ij}=1$  and $\sum^N_{j=1}[P(t)]_{ij}=1$ for all $t$ and $i,j\in\mathcal{V}$.
\end{ass}
\begin{ass}\label{ass2}
The time-varying network $(\mathcal{V},E(t),P(t))$ is $B$ connectivity. There exists a positive integer $B$ such that the graph
\begin{eqnarray*}
(\mathcal{V},\cup^{(c+1)B}_{t=cB+1}E(t))
\end{eqnarray*}
 is strongly connected for any $c\geq 0$. There is a $\theta\in(0,1)$ such that $[P(t)]_{ij}\geq \theta$ if $(i,j)\in E(t)$ and $[P(t)]_{jj}\geq \theta$ for all $j\in \mathcal{V}$.

 \end{ass}

Assumptions $\ref{ass1}$ and $\ref{ass2}$ are widely used in distributed optimization over time-varying communication network. In this paper, we define the transition matrix $P(t,s)=\prod^{0}_{i=s-t}P(s-i)$ for $t\geq s$ and $P(s,s+1)=I_n$  for any $s\geq 0$. The following Lemma \ref{Nedic11}  is critical in the analysis of time-varying communication network.

\begin{lem}\label{Nedic11}\cite{Nedic1}
If  Assumptions \ref{ass1} and \ref{ass2} are satisfied, then we have
\begin{eqnarray*}
|[P(m,n)]_{ij}-\frac{1}{N}|\leq \omega\gamma^{m-n},
\end{eqnarray*}
for all $i,j\in \mathcal{V}$ and $m,n$ satisfying $m\geq n\geq 1$, where $\omega=(1-\frac{\theta}{4N^2})^{-2}$ and $\gamma=(1-\frac{\theta}{4N^2})^{\frac{1}{B}}$.
\end{lem}


In this paper, we will  develop a mirror descent algorithm to solve the distributed optimization problem. We consider a $\sigma_{\phi}$ strongly convex and  continuously differentiable  distance generating  function $\phi(x):\mathbb{R} ^n\rightarrow \mathbb{R}$ and the Bregman divergence $V_{\phi}(a,b):\mathbb{R} ^n\times\mathbb{R} ^n\rightarrow \mathbb{R}$ associated with $\phi$ as
\begin{eqnarray}\label{BregmandivergenceV}
V_{\phi}(a,b)=\phi(a)-\phi(b)-\langle \nabla\phi(b),a-b \rangle.
\end{eqnarray}
We will use the following standard assumptions (e.g. \cite{Nedic5}) on the distance generating function $\phi(x)$.
\begin{ass} The gradient of the distance generating function $\phi(x)$ is Lipschitz continuity with constant $L_{\phi}$
\begin{eqnarray*}
||\nabla \phi(a)-\nabla\phi(b)||\leq L_{\phi}||a-b||,
\end{eqnarray*}
for any $a,b\in \mathbb{R} ^n$.
\end{ass}

The following separate convexity assumption on Bregman divergence is standard in the study of  distributed mirror descent algorithms (e.g. \cite{Yuan2, Yuan5}).
\begin{ass}
The Bregman divergence $V_{\phi}(x,y)$ satisfies the separate convexity. For any vector $a\in\mathbb{R}^n$ and a sequence vectors $\{b_i\}^N_{j=1}\in\mathbb{R}^n$, we have
\begin{eqnarray*}
V_{\phi}(a,\sum^N_{j=1}w_jb_j)\leq \sum^N_{i=1}w_jV_{\phi}(a,b_j),
\end{eqnarray*}
where $\sum^N_{j=1}w_j=1$.
\end{ass}

The following boundedness assumption is always used for (sub)gradient and Bregman divergence (e.g. \cite{Yuan5}).
\begin{ass}\label{assump1}
The norm of  $g_i(x)\in \partial f_i(x)$ and  $V_{\phi}(x,y)$ are uniformly bounded on $\mathcal{X}$ with
\begin{eqnarray}
\sup_{x\in \mathcal{X},1\leq i\leq N}||g_i(x)||&\leq& G,\label{Bound111}\\
\sup_{x\in\mathcal{X},y\in\mathcal{X}}V_{\phi}(x,y)&\leq& D_{\phi}.\label{Bound222}
\end{eqnarray}
\end{ass}

The distance generating function $\phi$ is $\sigma_{\phi}$-strongly convex and we know that
\begin{eqnarray}\label{distancegenerating}
\phi(x)\geq \phi(y)+\langle \nabla\phi(y),x-y \rangle+\frac{\sigma_{\phi}}{2}||x-y||^2.
\end{eqnarray}
For $\forall x,y\in \mathcal{X}$, take inequality (\ref{Bound222}), (\ref{distancegenerating}) into equation (\ref{BregmandivergenceV}), we have
\begin{eqnarray*}
\frac{\sigma_{\phi}}{2}||x-y||^2\leq V_{\phi}(x,y)\leq D_{\phi}
\end{eqnarray*}
and then  the bound of $\sup_{x,y\in\mathcal{X}}||x-y||$ is obtained
\begin{eqnarray}\label{boundofsupxy}
\sup_{x,y\in\mathcal{X}}||x-y||\leq\sqrt{\frac{2D_{\phi}}{\sigma_{\phi}}}.
\end{eqnarray}

\section{Mirror Descent Algorithm with Delayed Subgradient Information under Adaptive Quantization}

Due to the  limited communication capacity and delayed subgradient information,   the appropriate  quantizer needs to be designed.  However, the quantization error is caused after receiving the quantized information. In some previous works (e.g. \cite{Nedic4,Yuanquantization}), due to the  non-decreasing  quantization error's upper bounds, the convergence of the algorithms are inexact and the solutions are suboptimal.  In this paper, we design a mirror descent algorithm (proposed in \textbf{Algorithm 1}) with delayed subgradient information under adaptive quantization. The quantization error under adaptive quantizer decreases to $0$, which is different from that in \cite{Nedic4,Yuanquantization}. \textbf{Algorithm 1}'s optimal convergence rate  $O(1/\sqrt{T})$  can be obtained by selecting appropriate parameters.




Some notations in \textbf{Algorithm 1} need to be introduced. $x_i(t)$ and $z_i(t)$ are the state and quantization mid value of agent $i$ at time $t$, respectively; $d(t)$ is the quantization interval size for all agents at time $t$; The quantizer $\widehat{Q}$ is defined in equation (\ref{quantization}). $y_i(t)$ is the weighted sum of quantized information received from agent $i$'s all neighbors. The non-increasing positive sequences $\{\alpha(t)\}$ and $\{\beta(t)\}$ are defined as the sequences of stepsize and quantization parameter, respectively, where $0<\alpha(t), \beta(t)<1$ for any $t$ and $\alpha(t), \beta(t)$ are decreasing to $0$. Note that $\{\alpha(t)\}$ and $\{\beta(t)\}$ are pre-determined and need to satisfy some conditions, which will be discussed in detail in  Section \uppercase\expandafter{\romannumeral5}.

\begin{algorithm}[H]\label{alg-1}
\caption{Mirror Descent Algorithm with Delayed Subgradient Information under Adaptive Quantization}
  Initialize: $x_i(0)\in \mathcal{X}$, $z_i(0)=x_i(0)$, $\tilde{y}_i(-\tau)$, $\tilde{y}_i(-\tau+1)$, $\cdots$, $\tilde{y}_i(-1)\in \mathcal{X}$ with $i=1,\cdots,N$, stepsize sequence $\{\alpha(t)\}$ and quantization parameter sequence $\{\beta(t)\}$;\\
1: for $t=0,1,2,\cdots,$\\
2: $\quad$Compute $d(t)=\frac{G\alpha(t)\beta(t)}{\sigma_{\phi}}\textbf{1}$;\\
3: $\quad$for $i=1$:$i\leq N$:$i++$\\
4: $\quad\quad$Agent $i$ receives  quantized values $\widehat{Q}(z_j(t),d(t),x_j(t))$ from  all neighbors $j\in N_i(t)$;\\
5: $\quad\quad$Update $y_i(t)=\sum^N_{j=1}[P(t)]_{ij}\widehat{Q}(z_j(t),d(t),x_j(t))$;\\
6: $\quad\quad$Compute $\tilde{y}_i(t)=Proj_{x\in\mathcal{X}}\{y_i(t)\}$;\\
7: $\quad\quad$Compute the subgradient $g_i(t-\tau)\in \partial f_i(\tilde{y}_i(t-\tau))$;\\
8: $\quad\quad$Update $z_i(t+1)$ and $x_i(t+1)$ as follow:
 \begin{eqnarray*}
&&z_i(t+1)=\arg \min_{x\in\mathcal{X}}\{\langle g_i(t-\tau),x\rangle+\frac{V_{\phi}(x,\tilde{y}_i(t))}{\alpha(t+1)(1-\beta(t+1))}\}\\
&&x_i(t+1)=\arg \min_{x\in\mathcal{X}}\{\langle g_i(t-\tau),x\rangle+\frac{ V_{\phi}(x,\tilde{y}_i(t))}{\alpha(t+1)}\}
\end{eqnarray*}
9:$\quad$end for\\
10:end for
\end{algorithm}


The following key inequalities Lemmas \ref{Ghadimi} and \ref{fallinto} are provided to show that the state $x_i(t)$ falls into the quantization interval $[z_i(t)-d(t),z_i(t)+d(t)]$ at each iteration, which play a crucial role in the derivation of our main results.

\begin{lem}\label{Ghadimi}\cite{Ghadimi}
For the Bregman divergence $V_{\phi}(x,z)$ and
\begin{eqnarray*}
x^{+}_1&=&\arg\min_{x\in\mathcal{X}}\{\alpha\langle g_1,x\rangle+V_{\phi}(x,y)\},\\
x^{+}_2&=&\arg\min_{x\in\mathcal{X}}\{\alpha\langle g_2,x\rangle+V_{\phi}(x,y)\},
\end{eqnarray*}
where for any $\alpha>0$, $g_1,g_2\in R^n$ and $y\in \mathcal{X}$. Then we have
\begin{eqnarray*}
||x^{+}_2-x^{+}_1||\leq \frac{\alpha}{\sigma_{\phi}}||g_2-g_1||,
\end{eqnarray*}
where $\sigma_{\phi}$ is the modulus of strong convexity of distance generating function $\phi$.
\end{lem}

\begin{lem}\label{fallinto}
The $x_i(t)$ generated by \textbf{Algorithm 1} falls into the quantization intervals $[z_i(t)-d(t),z_i(t)+d(t)]$ for  $i=1,2,\cdots,N$ at each iteration.
\end{lem}
\begin{proof}
From Lemma \ref{Ghadimi}, for any $i=1,2,\cdots,N$ and $j=1,2,\cdots,n$, we know that
\begin{eqnarray*}
||x_i(t)-z_i(t)||&\leq&\frac{||g(t-1-\tau)||\alpha(t)\beta(t)}{\sigma_{\phi}}\\
&\leq&\frac{G\alpha(t)\beta(t)}{\sigma_{\phi}}\\
&=&[d(t)]_j.
\end{eqnarray*}
Hence, for any $j=1,2,\cdots,n$, we have
\begin{eqnarray*}%
|[x_i(t)-z_i(t)]_j|\leq||x_i(t)-z_i(t)||\leq [d(t)]_j,
\end{eqnarray*}
which is equivalent to
\begin{eqnarray*}
z_i(t)-d(t)\preceq x_i(t)\preceq z_i(t)+d(t).
\end{eqnarray*}
Therefore, $x_i(t)$ generated by \textbf{Algorithm 1} will fall into quantization intervals $[z_i(t)-d(t),z_i(t)+d(t)]$ for  $i=1,2,\cdots,N$ at each iteration.
\end{proof}

From Lemma \ref{fallinto} and inequality (\ref{quantizationbound}), we know that the quantization error satisfies
\begin{eqnarray}\label{quantizationerrorupperbound1}
||x_i(t)-\widehat{Q}(z_i(t),d(t),x_i(t))||\leq \sqrt{n}||d(t)||.
\end{eqnarray}

\noindent The quantization interval $d(t+1)$ in \textbf{Algorithm 1} satisfies
\begin{eqnarray}\label{quantizationerrorupperbound2}
||d(t)||=\frac{G\sqrt{n}}{\sigma_{\phi}}\alpha(t)\beta(t).
\end{eqnarray}
Both stepsize $\alpha(t)$ and quantization parameter $\beta(t)$ are to be designed to decrease to $0$ and then the quantization error $||x_i(t)-\widehat{Q}(z_i(t),d(t),x_i(t))||$ will decrease to $0$. 


\begin{rmk}
The technical Lemma \ref{Ghadimi} and Lemma \ref{fallinto} help us overcome the difficulty of designing the adaptive quantizer in non-Euclidean projection based distributed optimization algorithm. With Lemma \ref{Ghadimi}, we can  find the appropriate quantization mid value $z_i(t)$ and quantization interval size $d(t)$ to make the state $x_i(t)\in[z_i(t)-d(t),z_i(t)+d(t)]$. The above technical lemmas indicates the essential distinction with the existing works such as   \cite{Nedic4,Pu1,Yi,Doan}.
\end{rmk}
\section{Main Results}
In this section, we will analyze the convergence of distributed mirror descent algorithm with delayed subgradient information under adaptive quantization and discuss the convergence rate of different   stepsize $\alpha(t)$ and quantization parameter $\beta(t)$.
\subsection{Convergence Analysis}

The Steps $5$, $6$ and $8$ of \textbf{Algorithm 1} are shown as follow:
\begin{eqnarray}
\label{algorithm666111}\quad\quad y_i(t)\quad&=&\sum^N_{j=1}[P(t)]_{ij}\widehat{Q}(z_j(t),d(t),x_j(t)),\\
 \label{algorithm666222}\quad\quad \tilde{y}_i(t)\quad&=&Proj_{\mathcal{X}}\{y_i(t)\},\\
x_i(t+1)&=&\arg\min_{x\in \mathcal{X}}\{\langle g_i(t-\tau),x\rangle +\frac{V_{\phi}(x,\tilde{y}_i(t))}{\alpha(t+1)}\},\label{algorithm666333}
\end{eqnarray}
where $P[t]$ is time-varying double stochastic communication matrix at time $t$, $\tilde{y}_i(t)$ is  projection of $y_i(t)$ onto $\mathcal{X}$, $V_{\phi}$ is Bregman divergence and quantizer $\widehat{Q}$ is defined as equation (\ref{quantization}).

We will use the upper bound of quantization error to  analyze the convergence rate of \textbf{Algorithm 1}. In order to present clearly, we use $e_j(t)$ and $p_i(t)$ to denote quantization error  and projection error,  respectively with $i,j=1,2,\cdots, N$,
 \begin{eqnarray}
 \label{quantizationerror} e_j(t)&=&\widehat{Q}(z_j(t),d(t),x_j(t))-x_j(t),\\
 \label{projectionerror} p_i(t)&=&\tilde{y}_i(t)-y_i(t).
 \end{eqnarray}
 \noindent The equivalent forms of (\ref{algorithm666111}), (\ref{algorithm666222}) and (\ref{algorithm666333})   are
\begin{eqnarray}\label{equation-2}
\tilde{y}_i(t)&=&\sum^N_{j=1}[P(t)]_{ij}(x_j(t)+e_j(t))+p_i(t),\label{equation-2}\\
x_i(t+1)&=&\arg\min_{x\in \mathcal{X}}\{\langle g_i(t-\tau),x\rangle +\frac{V_{\phi}(x,\tilde{y}_i(t))}{\alpha(t+1)}\}\label{equation-3}.
\end{eqnarray}
From inequality (\ref{quantizationerrorupperbound1}) and equations (\ref{quantizationerrorupperbound2}), (\ref{quantizationerror}), hence we know that
\begin{eqnarray}\label{analysisofquantizationerror}
||e_j(t)||=||\widehat{Q}(z_j(t),d(t),x_j(t))-x_j(t)||\leq\frac{Gn}{\sigma_{\phi}}\alpha(t)\beta(t),
\end{eqnarray}
for $j=1,2,\cdots,N$.  We define that
\begin{eqnarray}\label{quantizationerrorupperbound3}
E(t)=\frac{Gn}{\sigma_{\phi}}\alpha(t)\beta(t)
\end{eqnarray}
and from (\ref{analysisofquantizationerror}) and (\ref{quantizationerrorupperbound3}), then we have
\begin{eqnarray}\label{qwerbound}
||e_j(t)||\leq E(t).
\end{eqnarray}
Similarly, the upper bound of projection error $||p_i(t)||$ should be obtained in order to analyze \textbf{Algorithm 1}'s convergence. We have the following Lemma \ref{lemma-4} to obtain the projection error's upper bound.
\begin{lem}\label{lemma-4}
 The Euclidean norm of projection error $p_i(t)$ satisfies
\begin{eqnarray}\label{Euclideannormofprojectionerror}
||p_i(t)||\leq 2NE(t),
\end{eqnarray}
where $E(t)$ is defined as (\ref{quantizationerrorupperbound3}).
\end{lem}

\begin{proof}
From equation (\ref{algorithm666111}), (\ref{quantizationerror}) and (\ref{projectionerror}), we have
\begin{eqnarray}
&&||p_i(t)||\nonumber\\
&=&||\tilde{y}_i(t)-y_i(t)||\nonumber\\
&=&||\tilde{y}_i(t)-\sum^N_{j=1}[P(t)]_{ij}x_j(t)-\sum^N_{j=1}[P(t)]_{ij}e_j(t)||\nonumber\\
&\leq&||\tilde{y}_i(t)-\sum^N_{j=1}[P(t)]_{ij}x_j(t)||+||\sum^N_{j=1}[P(t)]_{ij}e_j(t)||\label{projectionerrorinequality1}\\
&\leq&||y_i(t)-\sum^N_{j=1}[P(t)]_{ij}x_j(t)||+NE(t)\label{projectionerrorinequality2}\\
&=&2NE(t),\nonumber
\end{eqnarray}
where the  inequality (\ref{projectionerrorinequality1}) is obtained from triangle inequality and  inequality (\ref{projectionerrorinequality2}) is obtained from   projection theorem, respectively.
\end{proof}

We define  Bregman projection error as
\begin{eqnarray}\label{Bregmanprojectionerror}
\varepsilon_i(t)=x_i(t+1)-\tilde{y}_i(t)
\end{eqnarray}
and need to obtain the Bregman projection error's upper bound.  We have the following Lemma \ref{Bregman} about the upper bound of $||\varepsilon_i(t)||$.
\begin{lem}\label{Bregman}
Bregman projection error $\varepsilon_i(t)$ with $i=1,2,\cdots, N$ satisfies
\begin{eqnarray*}
||\varepsilon_i(t)||\leq \frac{G\alpha(t)}{\sigma_{\phi}}.
\end{eqnarray*}
\end{lem}
\begin{proof}
The first order optimality of $x_i(t+1)$ implies, for $\forall x\in\mathcal{X}$, we have
\begin{eqnarray*}
\langle \alpha(t)g_i(t-\tau)+\nabla \phi(x_i(t+1))-\nabla\phi(\tilde{y}_i(t)),x-x_i(t+1)\rangle\geq0.
\end{eqnarray*}
Substitute $\tilde{y}_i(t)\in\mathcal{X}$ into above inequality, and we have
\begin{eqnarray}\label{optimality2}
\langle \alpha(t)g_i(t-\tau)+\nabla \phi(x_i(t+1))-\nabla\phi(\tilde{y}_i(t)),\tilde{y}_i(t)-x_i(t+1)\rangle\geq0.
\end{eqnarray}
Rearrange the terms in (\ref{optimality2}) and we have
\begin{eqnarray}\label{strongconvex}
\langle \alpha(t)g_i(t-\tau),\tilde{y}_i(t)-x_i(t+1) \rangle&\geq& \langle\nabla \phi(x_i(t+1))-\nabla\phi(\tilde{y}_i(t)), x_i(t+1)-\tilde{y}_i(t)\rangle\nonumber\\
&\geq& \sigma_{\phi}||x_i(t+1)-\tilde{y}_i(t)||^2,
\end{eqnarray}
where the second inequality is derived from $\sigma_{\phi}$ strongly convex of $\phi(x)$. From Cauchy inequality, we have
\begin{eqnarray}\label{cauchy}
\alpha(t)||g_i(t-\tau)||||x_i(t+1)-\tilde{y}_i(t)||\geq \langle \alpha(t)g_i(t-\tau),\tilde{y}_i(t)-x_i(t+1) \rangle.
\end{eqnarray}
Combining (\ref{strongconvex}), (\ref{cauchy}) and subgradient's upper bound (\ref{Bound111}), we have
\begin{eqnarray}\label{howtogetbregman}
\sigma_{\phi}||x_i(t+1)-\tilde{y}_i(t)||^2\leq G\alpha(t)||x_i(t+1)-\tilde{y}_i(t)||.
\end{eqnarray}
Therefore, the upper bound of Bregman projection error $||\varepsilon_i(t)||$ is obtained from equation (\ref{Bregmanprojectionerror}) and inequality (\ref{howtogetbregman}).
\begin{eqnarray*}
||\varepsilon_i(t)||\leq\frac{G\alpha(t)}{\sigma_{\phi}}.
\end{eqnarray*}
The proof is completed.
\end{proof}

\begin{rmk}
The Step $6$ in \textbf{Algorithm 1} is necessary. Each agent receives the quantized value from its neighbors (shown in Step $5$). However, the quantized value $\widehat{Q}(z_j(t),d(t),x_j(t))$ of $x_j(t)$ may not belong to $\mathcal{X}$. Thus, the value $y_i(t)$ also may not belong to $\mathcal{X}$. In order to obtain the upper bound of Bregman projection error, we need to calculate the projection of $y_i(t)$ onto $\mathcal{X}$. Otherwise, inequality (\ref{optimality2}) is not true.
\end{rmk}

In order to derive the expression of $x_i(t)$ from equation (\ref{equation-2}), we define that $c_i(t)=\sum^N_{j=1}[P(t)]_{ij}(e_j(t)+p_i(t))$ and  have
\begin{eqnarray*}
x_i(t)&=&\sum^t_{s=1}\sum^N_{j=1}[P(t-1,s)]_{ij}(c_j(s-1)+\varepsilon_j(s-1))+\sum^N_{j=1}[P(t-1,0)]_{ij}x_j(0).
\end{eqnarray*}
The average state of all nodes at time $t$ is $\bar{x}(t)=\frac{1}{N}\sum^N_{i=1}x_i(t)$ and we have
\begin{eqnarray*}
\bar{x}(t)=\frac{1}{N}\sum^t_{s=1}\sum^N_{j=1}(c_j(s-1)+\varepsilon_j(s-1))+\frac{1}{N}\sum^N_{j=1}x_j(0).
\end{eqnarray*}
For each agent $l\in \mathcal{V}$, at iteration $T$, we consider  the following variable
\begin{eqnarray*}
\hat{x}_l(T)=\frac{1}{T}\sum^{T}_{t=1}x_l(t).
\end{eqnarray*}
We will analyze the property of $f(\hat{x}_l(T))-f(x^*)$ and derive its upper bound.  From the convexity of $f(x)$, we have
\begin{eqnarray}\label{veryveryImportant}
f(\hat{x}_l(T))-f(x^*)&\leq&\frac{1}{T}\sum^T_{t=1}f(x_l(t)) -f(x^*)\nonumber\\
&=&\frac{1}{T}\sum^{T+\tau}_{t=1+\tau}f(x_l(t-\tau)) -f(x^*).
\end{eqnarray}
Thus, we need to derive the upper bound of $\frac{1}{T}\sum^{T+\tau}_{t=1+\tau}f(x_l(t-\tau)) -f(x^*)$.
\begin{eqnarray}\label{cccccc}
\frac{1}{N}\sum^N_{i=1}\langle g_i(t-\tau),\tilde{y}_i(t-\tau)-x^*\rangle&\geq& \frac{1}{N}\sum^N_{i=1}\{f_i(\tilde{y}_i(t-\tau))-f_i(x^*)\}\nonumber\\
&=&\frac{1}{N}\sum^N_{i=1}f_i(\tilde{y}_i(t-\tau))-f(x^*).
\end{eqnarray}
The inequality in (\ref{cccccc}) is obtained from definition of subgradient $g_i(t-\tau)\in \partial f_i(\tilde{y}_i(t-\tau))$. For  $\frac{1}{N}\sum^N_{i=1}f_i(\tilde{y}_i(t-\tau))$ and any $1\leq l\leq N$, we have
\begin{eqnarray}
&&\frac{1}{N}\sum^N_{i=1}f_i(\tilde{y}_i(t-\tau))\nonumber\\
&=&\frac{1}{N}\sum^N_{i=1}\Big\{f_i(x_l(t-\tau))+f_i(\tilde{y}_i(t-\tau))-f_i(x_l(t-\tau))\Big\}\nonumber\\
&\geq& \frac{1}{N}\sum^N_{i=1}\langle g_i(x_l(t-\tau)), \tilde{y}_i(t-\tau)-x_l(t-\tau)\rangle +f(x_l(t-\tau))\label{abcabcabcd}\\
&\geq& f(x_l(t-\tau))-\frac{G}{N}\sum^N_{i=1}||\tilde{y}_i(t-\tau)-x_l(t-\tau)||,\label{aaaaaa}
\end{eqnarray}
the  inequality (\ref{abcabcabcd}) is obtained from definition of $g_i(x_l(t-\tau))\in\partial f_i(x_l(t-\tau))$ and the inequality (\ref{aaaaaa}) is obtained from Cauchy inequality and $||g_i(x_l(t-\tau))||\leq G$. From equation (\ref{equation-2}), we have
\begin{eqnarray}
&&||\tilde{y}_i(t-\tau)-x_l(t-\tau)||\nonumber\\
&=&||\sum^N_{j=1}[P(t-\tau)]_{ij}(x_j(t-\tau)-x_l(t-\tau))+\sum^N_{j=1}[P(t-\tau)]_{ij}e_j(t-\tau)+p_i(t-\tau)||\nonumber\\
&\leq&\sum^N_{j=1}[P(t-\tau)]_{ij}||x_j(t-\tau)-x_l(t-\tau)||+\sum^N_{j=1}[P(t-\tau)]_{ij}||e_j(t-\tau)||+||p_i(t-\tau)||\label{bbbbbbbbb}\\
&\leq&\sum^N_{j=1}[P(t-\tau)]_{ij}||x_j(t-\tau)-x_l(t-\tau)||+3NE(t-\tau),\label{bbbbb}
\end{eqnarray}
where inequality (\ref{bbbbbbbbb}) is obtained from triangle inequality and inequality (\ref{bbbbb}) is obtained from inequality (\ref{qwerbound}) and (\ref{Euclideannormofprojectionerror}). Therefore, from inequalities (\ref{aaaaaa}) and (\ref{bbbbb}), we have
\begin{eqnarray}\label{dddddd}
&&\frac{1}{N}\sum^N_{i=1}f_i(\tilde{y}_i(t-\tau))\nonumber\\
&\geq& f(x_l(t-\tau))-\frac{G}{N}\sum^N_{i=1}||\tilde{y}_i(t-\tau)-x_l(t-\tau)||\nonumber\\
&\geq&  f(x_l(t-\tau))-\frac{G}{N}\sum^N_{i=1}3NE(t-\tau)-\frac{G}{N}\sum^N_{i=1}\sum^N_{j=1}[P(t-\tau)]_{ij}||x_j(t-\tau)-x_l(t-\tau)||\nonumber\\
&=&f(x_l(t-\tau))-3NGE(t-\tau)-\frac{G}{N}\sum^N_{j=1}||x_j(t-\tau)-x_l(t-\tau)||.
\end{eqnarray}
Substitute (\ref{dddddd}) into (\ref{cccccc}) and we have
\begin{eqnarray}\label{eeeeee}
&&\frac{1}{N}\sum^N_{i=1}\langle g_i(t-\tau),\tilde{y}_i(t-\tau)-x^*\rangle\nonumber\\
&\geq& f(x_l(t-\tau)) -f(x^*)-3NGE(t-\tau)-\frac{G}{N}\sum^N_{j=1}||x_j(t-\tau)-x_l(t-\tau)||.
\end{eqnarray}
After summing up inequality (\ref{eeeeee}) from $t=\tau+1$ to $T+\tau$ and  dividing both side by $T$,   we have
\begin{eqnarray}\label{ffffff}
&&\frac{1}{NT}\sum^{T+\tau}_{t=1+\tau}\sum^N_{i=1}\langle g_i(t-\tau),\tilde{y}_i(t-\tau)-x^*\rangle\nonumber\\
&\geq& \frac{1}{T}\sum^{T+\tau}_{t=1+\tau}f(x_l(t-\tau)) -f(x^*)-\frac{3NG}{T}\sum^{T+\tau}_{t=1+\tau}E(t-\tau)\nonumber\\
&&-\frac{G}{TN}\sum^{T+\tau}_{t=1+\tau}\sum^N_{j=1}||x_j(t-\tau)-x_l(t-\tau)||.
\end{eqnarray}
Rearrange terms in inequality (\ref{ffffff}) and we can derive the upper bound of the following terms in (\ref{ffffff}),
\begin{eqnarray}\label{gggggg}
&&\frac{1}{T}\sum^{T+\tau}_{t=1+\tau}f(x_l(t-\tau)) -f(x^*)\nonumber\\
&\leq& \frac{1}{NT}\sum^{T+\tau}_{t=1+\tau}\sum^N_{i=1}\langle g_i(t-\tau),\tilde{y}_i(t-\tau)-x^*\rangle+\frac{G}{TN}\sum^{T+\tau}_{t=1+\tau}\sum^N_{j=1}||x_j(t-\tau)-x_l(t-\tau)||\nonumber\\
&&+\frac{3NG}{T}\sum^{T+\tau}_{t=1+\tau}E(t-\tau).
\end{eqnarray}
Substitute inequality (\ref{gggggg}) into inequality (\ref{veryveryImportant}) and we have
\begin{eqnarray}\label{Important}
&&f(\hat{x}_l(T))-f(x^*)\nonumber\\
&\leq&\frac{1}{NT}\sum^{T+\tau}_{t=1+\tau}\sum^N_{i=1}\langle g_i(t-\tau),\tilde{y}_i(t-\tau)-x^*\rangle+\frac{3NG}{T}\sum^{T+\tau}_{t=1+\tau}E(t-\tau)\nonumber\\
&&+\frac{G}{TN}\sum^{T+\tau}_{t=1+\tau}\sum^N_{j=1}||x_j(t-\tau)-x_l(t-\tau)||\nonumber\\
&=&\frac{1}{NT}\sum^{T+\tau}_{t=1+\tau}\sum^N_{i=1}\langle g_i(t-\tau),\tilde{y}_i(t)-x^*\rangle+\frac{1}{NT}\sum^{T+\tau}_{t=1+\tau}\sum^N_{i=1}\langle g_i(t-\tau),\tilde{y}_i(t-\tau)-\tilde{y}_i(t)\rangle\nonumber\\
&&+\frac{G}{TN}\sum^{T+\tau}_{t=1+\tau}\sum^N_{j=1}||x_j(t-\tau)-x_l(t-\tau)||+\frac{3NG}{T}\sum^{T+\tau}_{t=1+\tau}E(t-\tau).
\end{eqnarray}
In order to obtain the upper bound of $f(\hat{x}_l(T))-f(x^*)$, we need   Lemmas \ref{lemma-6}-\ref{lemma-2} to estimate the upper bound of the first, second and third term in (\ref{Important}), respectively.

\begin{lem}\label{lemma-6}
For the first term in (\ref{Important}), we have
\begin{eqnarray}\label{AA}
\frac{1}{NT}\sum^{T+\tau}_{t=1+\tau}\sum^N_{i=1}\langle g_i(t-\tau),\tilde{y}_i(t)-x^*\rangle\leq \mathcal{B}_1.
\end{eqnarray}
where $\mathcal{B}_1=\sqrt{\frac{2D_{\phi}}{\sigma_{\phi}}}\frac{(2N+1)L_{\phi}}{T}\sum^{T+\tau}_{t=1+\tau}\frac{E(t)}{\alpha(t)}+\frac{G^2}{2\sigma_{\phi}T}\sum^{T+\tau}_{t=1+\tau}\alpha(t)+\frac{D_{\phi}}{T\alpha(T+\tau)}+L_{\phi}(8N^2+2)\frac{1}{T}\sum^{T+\tau}_{t=1+\tau}\frac{E^2(t)}{\alpha(t)}.$
\end{lem}
\begin{proof}
See Appendix A.
\end{proof}

\begin{lem}\label{lemma-7}
For the second term in (\ref{Important}), we have
\begin{eqnarray}\label{BB}
\frac{1}{NT}\sum^{T+\tau}_{t=1+\tau}\sum^N_{i=1}\langle g_i(t-\tau),\tilde{y}_i(t-\tau)-\tilde{y}_i(t)\rangle\leq \mathcal{B}_2.
\end{eqnarray}
where  $\mathcal{B}_2=\frac{(2B+1)G^2\tau}{\sigma_{\phi}T}\sum^{T+\tau}_{t=0}\alpha(t)+\frac{(6B+3)N G\tau}{T}\sum^{T+\tau}_{t=0}E(t)+\frac{2N G\tau\omega A}{1-\gamma}\frac{1}{T}$.
\end{lem}

\begin{proof}
See Appendix B.
\end{proof}

\begin{lem}\label{lemma-2}
We define that $A=\sum^N_{j=1}||x_j(0)||$, $B=(2N+\frac{N^2\omega}{1-\gamma})$ and  have
\begin{eqnarray*}
\sum^T_{t=1}\sum^N_{i=1}||x_i(t)-\bar{x}(t)||\leq \frac{N\omega}{1-\gamma}A+B\sum^T_{t=0}\{\frac{G\alpha(t)}{\sigma_{\phi}}+3NE(t)\}
\end{eqnarray*}
\begin{eqnarray*}
\sum^T_{t=1}\sum^N_{i=1}||x_i(t)-x_j(t)||\leq \frac{2N\omega}{1-\gamma}A+2B\sum^T_{t=0}\{\frac{G\alpha(t)}{\sigma_{\phi}}+3NE(t)\}
\end{eqnarray*}
\end{lem}
\begin{proof}
See Appendix C.
\end{proof}
From Lemma \ref{lemma-2},  for the third term in (\ref{Important}), we define that $\mathcal{B}_3=\frac{2GA\omega}{T(1-\gamma)}+\frac{6BG}{T}\sum^T_{t=0}E(t)+\frac{2B G^2}{\sigma_{\phi}NT}\sum^T_{t=0}\alpha(t)$ and  have
\begin{eqnarray}\label{CC}
\frac{G}{TN}\sum^{T+\tau}_{t=1+\tau}\sum^N_{j=1}||x_j(t-\tau)-x_l(t-\tau)||\leq \mathcal{B}_3.
\end{eqnarray}
With the above Lemmas \ref{lemma-6}-\ref{lemma-2}, we have established the upper bound of the first, second and third terms in (\ref{Important}) and then we have the following theorem.

\begin{thm}
Suppose that Assumptions $1$-$5$ holds and let the sequence $\{x_l(t)\}_{k\geq 1}$ for all $l\in \mathcal{V}$ be generated by \textbf{Algorithm 1}. Let $f(x^*)$ be the optimal value of problem (\ref{goal}). Then for any $l\in \mathcal{X}$, we have

\begin{eqnarray}\label{thm}
&&f(\hat{x}_l(T))-f(x^*)\nonumber\\
&\leq&\{\frac{(2B+1)G^2\tau}{\sigma_{\phi}}+\frac{G^2}{2\sigma_{\phi}}+\frac{2BG^2}{{N\sigma_{\phi}}}\}\frac{1}{T}\sum^{T+\tau}_{t=0}\alpha(t)+\{(6B+3)N G\tau+6BG+3NG\}\frac{1}{T}\sum^{T+\tau}_{t=0}E(t)\nonumber\\
&&+\{\frac{2NAG\tau\omega}{1-\gamma}+\frac{2AG\omega}{1-\gamma}\}\frac{1}{T}+\frac{L_{\phi}(8N^2+2)}{T}\sum^{T+\tau}_{t=0}\frac{E^2(t)}{\alpha(t)}+\frac{(2N+1)L_{\phi}}{T}\sqrt{\frac{2D_{\phi}}{\sigma_{\phi}}}\sum^{T+\tau}_{t=0}\frac{E(t)}{\alpha(t)}\nonumber\\
&&+\frac{D_{\phi}}{T\alpha(T+\tau)},
\end{eqnarray}
where $A=\sum^N_{j=1}||x_j(0)||$, $B=(2N+\frac{N^2\omega}{1-\gamma})$.
\end{thm}

\begin{proof}
Substitute inequalities (\ref{AA}), (\ref{BB}) and (\ref{CC}) into inequality (\ref{Important}),  and then we have
\begin{eqnarray*}
&&f(\hat{x}_l(T))-f(x^*)\nonumber\\
&\leq&\mathcal{B}_1+\mathcal{B}_2+\mathcal{B}_3+\frac{3NG}{T}\sum^{T+\tau}_{t=1+\tau}E(t-\tau)\nonumber\\
&\leq&\frac{(2B+1)G^2\tau}{\sigma_{\phi}}\frac{1}{T}\sum^{T+\tau}_{t=0}\alpha(t)+\frac{(6B+3)N G\tau}{T}\sum^{T+\tau}_{t=0}E(t)+\frac{1}{T}\frac{2N G\tau\omega}{1-\gamma}A+6BG\frac{1}{T}\sum^{T+\tau}_{t=0}E(t)\nonumber\\
&&+\frac{G^2}{2\sigma_{\phi}}\frac{1}{T}\sum^{T+\tau}_{t=0}\alpha(t)+\frac{D_{\phi}}{T\alpha(T+\tau)}+\frac{1}{T}\frac{2AG\omega}{1-\gamma}+L_{\phi}(8N^2+2)\frac{1}{T}\sum^{T+\tau}_{t=0}\frac{E^2(t)}{\alpha(t)}+\frac{3NG}{T}\sum^{T+\tau}_{t=0}E(t)\nonumber\\
&&+\frac{2BG^2}{{N\sigma_{\phi}}}\frac{1}{T}\sum^{T+\tau}_{t=0}\alpha(t)+(2N+1)L_{\phi}\sqrt{\frac{2D_{\phi}}{\sigma_{\phi}}}\frac{1}{T}\sum^{T+\tau}_{t=0}\frac{E(t)}{\alpha(t)},\nonumber
\end{eqnarray*}
Rearrange the terms and we can get inequality (\ref{thm}). Therefore, the upper bound of $f(\hat{x}_l(T))-f(x^*)$ is obtained. The proof is completed.
\end{proof}
\begin{rmk}
We get the upper bound of $f(\hat{x}_l(T))-f(x^*)$ in Theorem 1.  From the right hand side of inequality (\ref{thm}), we know that the properties of stepsize $\alpha(t)$ and quantization parameter $\beta(t)$  decide the convergence of \textbf{Algorithm 1}. Thus,  the appropriate conditions of $\alpha(t)$ and $\beta(t)$  need to be discussed. Next, we will discuss how to design stepsize and quantization parameter.
\end{rmk}

\subsection{How to Design Stepsize and Quantization Parameter}
In this subsection, we will show some conditions of stepsize $\alpha(t)$ and quantization parameter $\beta(t)$ for \textbf{Algorithm 1}'s convergence.
\begin{thm}\label{thm2}
When the stepsize $\alpha(t)$ and quantization parameter $\beta(t)$ satisfy the following conditions (\ref{AAAA2}), (\ref{BBBB2}) and (\ref{CCCCdsfsad2}),
\begin{eqnarray}
\lim_{T\rightarrow \infty}\frac{1}{T\alpha(T+\tau)}&=&0,\label{AAAA2}\\
\lim_{T\rightarrow \infty}\frac{1}{T}\sum^{T+\tau}_{t=1}\alpha(t)&=&0,\label{BBBB2}\\
\lim_{T\rightarrow \infty}\frac{1}{T}\sum^{T+\tau}_{t=1}\beta(t)&=&0,\label{CCCCdsfsad2}
\end{eqnarray}
then the convergence of the \textbf{Algorithm 1} is guaranteed.

\begin{proof}
Since $E(t)=\frac{Gn}{\sigma_{\phi}}\alpha(t)\beta(t)$ from equation (\ref{quantizationerrorupperbound3}) and we have
\begin{eqnarray}
\label{QQ1}\frac{1}{T}\sum^{T+\tau}_{t=0}E(t)&=&\frac{Gn}{\sigma_{\phi}}\frac{1}{T}\sum^{T+\tau}_{t=0}\alpha(t)\beta(t),\\
\label{QQ2}\frac{1}{T}\sum^{T+\tau}_{t=0}\frac{E^2(t)}{\alpha(t)}&=&\frac{G^2n^2}{\sigma^2_{\phi}}\frac{1}{T}\sum^{T+\tau}_{t=0}\alpha(t)\beta^2(t),\\
\label{QQ3}\frac{1}{T}\sum^{T+\tau}_{t=0}\frac{E(t)}{\alpha(t)}&=&\frac{Gn}{\sigma_{\phi}}\frac{1}{T}\sum^{T+\tau}_{t=0}\beta(t).
\end{eqnarray}

For the upper bound of $f(\hat{x}_l(T))-f(x^*)$,   substitute equations (\ref{QQ1}), (\ref{QQ2}) and (\ref{QQ3}) into inequality (\ref{thm}) and we have
\begin{eqnarray}\label{quantizeranalysis}
&&f(\hat{x}_l(T))-f(x^*)\nonumber\\
&\leq&\{\frac{(2B+1)G^2\tau}{\sigma_{\phi}}+\frac{G^2}{2\sigma_{\phi}}+\frac{2BG^2}{{N\sigma_{\phi}}}\}\frac{1}{T}\sum^{T+\tau}_{t=0}\alpha(t)+(2N+1)L_{\phi}\sqrt{\frac{2D_{\phi}}{\sigma_{\phi}}}\frac{Gn}{\sigma_{\phi}}\frac{1}{T}\sum^{T+\tau}_{t=0}\beta(t)\nonumber\\
&&+\{\frac{2NAG\tau\omega}{1-\gamma}+\frac{2AG\omega}{1-\gamma}\}\frac{1}{T}+\frac{D_{\phi}}{T\alpha(T+\tau)}+L_{\phi}(8N^2+2)\frac{G^2n^2}{\sigma^2_{\phi}}\frac{1}{T}\sum^{T+\tau}_{t=0}\alpha(t)\beta^2(t)\nonumber\\
&&+\{(6B+3)N G\tau+6BG+3NG\}\frac{Gn}{\sigma_{\phi}}\frac{1}{T}\sum^{T+\tau}_{t=0}\alpha(t)\beta(t).
\end{eqnarray}
From equations (\ref{BBBB2}) and (\ref{CCCCdsfsad2}), when $T$ is sufficiently large, we have
\begin{eqnarray*}
\frac{1}{T}\sum^{T+\tau}_{t=0}\alpha(t)\beta(t)&\leq&\frac{1}{T}\sum^{T+\tau}_{t=0}\alpha(t)\\
\frac{1}{T}\sum^{T+\tau}_{t=0}\alpha(t)\beta^2(t)&\leq&\frac{1}{T}\sum^{T+\tau}_{t=0}\alpha(t).
\end{eqnarray*}
Therefore, we have
\begin{eqnarray}
\lim_{T\rightarrow \infty}\frac{1}{T}\sum^{T+\tau}_{t=0}E(t)&=&\lim_{T\rightarrow \infty}\frac{Gn}{\sigma_{\phi}}\frac{1}{T}\sum^{T+\tau}_{t=0}\alpha(t)\beta(t)=0,\label{thmmm1}\\
\lim_{T\rightarrow \infty}\frac{1}{T}\sum^{T+\tau}_{t=0}\frac{E^2(t)}{\alpha(t)}&=&\lim_{T\rightarrow \infty}\frac{G^2n^2}{\sigma^2_{\phi}}\frac{1}{T}\sum^{T+\tau}_{t=0}\alpha(t)\beta^2(t)=0,\label{thmmm2}\\
\lim_{T\rightarrow \infty}\frac{1}{T}\sum^{T+\tau}_{t=0}\frac{E(t)}{\alpha(t)}&=&\lim_{T\rightarrow \infty}\frac{Gn}{\sigma_{\phi}}\frac{1}{T}\sum^{T+\tau}_{t=0}\beta(t)=0.\label{thmmm3}
\end{eqnarray}

From equations (\ref{AAAA2}), (\ref{BBBB2}), (\ref{thmmm1}), (\ref{thmmm2}), (\ref{thmmm3}), we know that each term in inequality (\ref{quantizeranalysis}) decreases to zero. Hence, we  have
\begin{eqnarray*}
\lim_{T\rightarrow \infty}f(\hat{x}_l(T))=f(x^*).
\end{eqnarray*}
The proof is completed.
\end{proof}
\end{thm}

From  Theorem 2, we know that the convergence of the \textbf{Algorithm 1} is guaranteed if stepsize $\alpha(t)$ and quantization parameter $\beta(t)$ satisfy the conditions (\ref{AAAA2}), (\ref{BBBB2}) and (\ref{CCCCdsfsad2}). The following Corollary \ref{cor6666} shows the equivalent condition of $\alpha(t)$ and  $\beta(t)$ to guarantee the convergence of \textbf{Algorithm 1}.
\begin{cor}\label{cor6666}
The equivalent forms of the conditions (\ref{AAAA2}), (\ref{BBBB2}) and (\ref{CCCCdsfsad2}) in Theorem \ref{thm2} are
\begin{eqnarray}
\lim_{T\rightarrow \infty}\frac{1}{T\alpha(T)}&=&0,\label{DDDD2}\\
\lim_{T\rightarrow \infty}\frac{1}{T}\sum^{T}_{t=1}\alpha(t)&=&0,\label{EEEE2}\\
\lim_{T\rightarrow \infty}\frac{1}{T}\sum^{T}_{t=1}\beta(t)&=&0.\label{FFFF2}
\end{eqnarray}


\end{cor}
\begin{proof}
See Appendix D.
\end{proof}

\begin{rmk}
Corollary \ref{cor6666} shows that  the convergence of \textbf{Algorithm 1} is guaranteed if the conditions (\ref{DDDD2}), (\ref{EEEE2}) and (\ref{FFFF2}) are satisfied. Furthermore, the time delay does not affect the conditions (\ref{DDDD2}), (\ref{EEEE2}) and (\ref{FFFF2}). This fact indicates that the time delay does not affect the convergence of \textbf{Algorithm 1}.
\end{rmk}

Furthermore, different  stepsize $\alpha(t)$ and quantization parameter $\beta(t)$ satisfying conditions (\ref{DDDD2}), (\ref{EEEE2}) and (\ref{FFFF2}) will affect the convergence rate of \textbf{Algorithm 1}. We have the following Corollary \ref{cor2}.

\begin{cor}\label{cor2}
 When stepsize $\alpha(t)$ and quantization parameter $\beta(t)$ satisfy the condition (\ref{AAAA2}), (\ref{BBBB2}) and (\ref{CCCCdsfsad2}), the  convergence rate is determined by $\max\{\frac{1}{T}\sum^{T+\tau}_{t=1}\beta(t),\frac{1}{T}\sum^{T+\tau}_{t=1}\alpha(t),\frac{1}{T\alpha({T+\tau})}\}$.
\end{cor}
\begin{proof}
When $T$ is sufficiently large, we have
\begin{eqnarray*}
\frac{1}{T}\sum^{T+\tau}_{t=0}\alpha(t)\beta(t)&\leq&\frac{1}{T}\sum^{T+\tau}_{t=0}\alpha(t),\\
\frac{1}{T}\sum^{T+\tau}_{t=0}\alpha(t)\beta^2(t)&\leq&\frac{1}{T}\sum^{T+\tau}_{t=0}\alpha(t).
\end{eqnarray*}
From condition (\ref{BBBB2}), we know that
\begin{eqnarray*}
\lim_{T\rightarrow \infty}\alpha(t)=0.
\end{eqnarray*}
Therefore, when $T$ is sufficiently large, we  have
\begin{eqnarray*}
\frac{1}{T}\leq \frac{1}{\alpha(T+\tau)T}.
\end{eqnarray*}
From Theorem \ref{thm2}, we know that the convergence rate of \textbf{Algorithm 1} is determined by
\begin{eqnarray}\label{rate11}
\max\{\frac{1}{T}\sum^{T+\tau}_{t=1}\beta(t),\frac{1}{T}\sum^{T+\tau}_{t=1}\alpha(t),\frac{1}{T\alpha({T+\tau})}\}.
\end{eqnarray}
\end{proof}

Corollary \ref{cor2} analyze the convergence rate of \textbf{Algorithm 1} under limited communication capacity. Next we will discuss the convergence rate of  \textbf{Algorithm 1} with stepsize $\alpha(t)=\frac{1}{(t+1)^{\rho_1}}$, quantization parameter $\beta(t)=\frac{1}{(1+t)^{\rho_2}}$, where $0<\rho_1, \rho_2<1$ .

\begin{cor}\label{cor22222222}
For the stepsize $\alpha(t)=\frac{1}{(t+1)^{\rho_1}}$, quantization parameter $\beta(t)=\frac{1}{(1+t)^{\rho_2}}$, where $0<\rho_1, \rho_2<1$, the convergence rate of \textbf{Algorithm 1} is $O(1/T^{\rho})$, where $\rho=\min\{\rho_1,\rho_2,1-\rho_1\}$.
\end{cor}
\begin{proof}
From Corollary \ref{cor2}, we need to analyze three terms $\frac{1}{T}\sum^{T+\tau}_{t=0}\alpha(t)$, $\frac{1}{T}\sum^{T+\tau}_{t=0}\beta(t)$ and $\frac{1}{T\alpha(T+\tau)}$ in (\ref{rate11}).
For the term $\frac{1}{T}\sum^{T+\tau}_{t=0}\alpha(t)$, we have
\begin{eqnarray}\label{cor4111}
\frac{1}{T}\sum^{T+\tau}_{t=0}\alpha(t)&\leq& \frac{1}{T}\{1+\int^{T+\tau+1}_{1}\frac{1}{(1+x)^{\rho_1}}dx\}\nonumber\\
&\leq&\frac{2^{1-\rho_1}}{1-\rho_1}\frac{1}{T^{\rho_1}}-\frac{\rho_1}{1-\rho_1}\frac{1}{T}.
\end{eqnarray}
For the term $\frac{1}{T}\sum^{T+\tau}_{t=0}\beta(t)$, similar to (\ref{cor4111}), we have
\begin{eqnarray}\label{cor4222}
\frac{1}{T}\sum^{T+\tau}_{t=0}\beta(t)\leq \frac{2^{1-\rho_2}}{1-\rho_2}\frac{1}{T^{\rho_2}}-\frac{\rho_2}{1-\rho_2}\frac{1}{T}.
\end{eqnarray}
For the term $\frac{1}{T\alpha(T+\tau)}$ and $T>\tau+1$, we have
\begin{eqnarray}\label{cor4333}
\frac{1}{T\alpha(T+\tau)}=\frac{(T+\tau+1)^{\rho_1}}{T}=\frac{2^{\rho_1}}{T^{1-\rho_1}}
\end{eqnarray}
We define that
\begin{eqnarray*}
\rho=\min\{\rho_1,\rho_2,1-\rho_1\},
\end{eqnarray*}
Thus the convergence rate is $O(1/T^{\rho})$. The proof is completed.
\end{proof}

The following table shows the convergence rate of different stepsize $\alpha(t)$ and quantization parameter $\beta(t)$.
\begin{table}[H]
\centering
\begin{tabular}{c|c|c|c}
\toprule
\diagbox{$\alpha(t)$}{$\beta(t)$}&$1/T^{0.1}$&$1/T^{0.5}$&$1/T^{0.9}$\\
\hline
$1/T^{0.25}$&$O(1/T^{0.1})$&$O(1/T^{0.25})$&$O(1/T^{0.25})$\\
\hline
$1/T^{0.5}$&$O(1/T^{0.1})$&$O(1/T^{0.5})$&$O(1/T^{0.5})$\\
\hline
$1/T^{0.75}$&$O(1/T^{0.1})$&$O(1/T^{0.25})$&$O(1/T^{0.25})$\\
\bottomrule
\end{tabular}
\end{table}

$\\$ Next we will discuss how to select appropriate stepsize $\alpha(t)$ and quantization parameter $\beta(t)$ to obtain the optimal convergence rate.
\begin{cor}
The optimal convergence rate is $O(1/\sqrt{T})$ as we select $\alpha(t)=\frac{1}{(t+1)^{\rho_1}}$ and $\beta(t)=\frac{1}{(t+1)^{\rho_2}}$, where $\rho_1=\frac{1}{2}$ and $\frac{1}{2}\leq\rho_2<1$.
\end{cor}


\begin{proof}
From Corollary \ref{cor22222222}, we know that the optimal convergence rate is $O(1/T^{\rho})$ and
\begin{eqnarray*}
\rho=\min\{\rho_1,1-\rho_1, \rho_2\}\leq \min\{\frac{1}{2},\rho_2\}\leq \frac{1}{2},
\end{eqnarray*}
where $\rho_1=\frac{1}{2}$ and then we have $\min\{\rho_1,1-\rho_1,\rho_2\}=\frac{1}{2}$. The optimal convergence rate is $O(1/\sqrt{T})$. The proof is completed.
\end{proof}

\begin{rmk}
As stated in the introduction, our algorithm's convergence rate is faster than that in \cite{Doan}. Also, the communication network in \cite{Doan} is static while we consider a class of time-varying communication network, which is more realistic. Furthermore, the Bregman divergence is considered instead of the classical Euclidean projection procedure in  \cite{Doan}, making the proposed algorithm more flexible when handling distributed optimization problems.
\end{rmk}

\begin{rmk}
The work in \cite{Li} investigates the mirror descent algorithm with delayed gradient and the convergence rate is $O(1/\sqrt{T})$. Compared with  \cite{Li}, we design an adaptive quantization method  and the convergence rate is $O(1/\sqrt{T})$. Furthermore, the assumptions in this paper are  easier to satisfy. In this paper, we assume that  objective function's subgradient has the upper bound, in contrast to this work, the objective function's gradient  in \cite{Li} is assumed to satisfy Lipschitz continuous, which is a stronger condition.
\end{rmk}


For a special case, when the communication between each agent is perfect, that is to say $E(t)=0$, we have the following Corollary \ref{cor3}.
\begin{cor}\label{cor3}
 The different stepsize $\alpha(t)=\frac{1}{(t+1)^{\rho_1}}$  with $0$$<$$\rho_1$$<$$1$ will affect the convergence rate under perfect communication. The convergence rate is determined by $\max\{\frac{1}{T}\sum^{T+\tau}_{t=1}\alpha(t),\frac{1}{T\alpha({T+\tau})}\}$. Furthermore, the optimal convergence rate is $O(1/\sqrt{T})$ and the corresponding stepsize is $\alpha(t)=\frac{1}{\sqrt{t+1}}$.
\end{cor}
\begin{proof}
When the communication between each agent is perfect, that is to say there is no quantization between communication. Then we have $E(t)=0$ and substitute it into inequality (\ref{thm}). That is
\begin{eqnarray*}
&&f(\hat{x}_l(T))-f(x^*)\nonumber\\
&\leq&\{\frac{(2B+1)G^2\tau}{\sigma_{\phi}}+\frac{G^2}{2\sigma_{\phi}}+\frac{2BG^2}{{N\sigma_{\phi}}}\}\frac{1}{T}\sum^{T+\tau}_{t=0}\alpha(t)+\{\frac{2NAG\tau\omega}{1-\gamma}+\frac{2AG\omega}{1-\gamma}\}\frac{1}{T}+\frac{D_{\phi}}{T\alpha(T+\tau)}.
\end{eqnarray*}
When $T$ is large enough, we  have
\begin{eqnarray*}
\frac{1}{T}\leq \frac{1}{\alpha(T+\tau)T}.
\end{eqnarray*}
Therefore, the convergence rate is determined by
\begin{eqnarray*}
\max\{\frac{1}{T}\sum^{T+\tau}_{t=1}\alpha(t),\frac{1}{T\alpha({T+\tau})}\}.
\end{eqnarray*}
From inequalities (\ref{cor4111}) and (\ref{cor4333}), we know that

\begin{eqnarray*}
\frac{1}{T}\sum^{T+\tau}_{t=0}\alpha(t)&\leq& \frac{2^{1-\rho_1}}{1-\rho_1}\frac{1}{T^{\rho_1}}-\frac{\rho_1}{1-\rho_1}\frac{1}{T},\\
\frac{1}{T\alpha(T+\tau)}&\leq&\frac{2^{\rho_1}}{T^{1-\rho_1}}.
\end{eqnarray*}
Therefore, the convergence rate is $O(1/T^{\rho})$, where $\rho=\min\{\rho_1,1-\rho_1\}$. We know that
\begin{eqnarray*}
\min\{\rho_1,1-\rho_1\}\leq \frac{1}{2}.
\end{eqnarray*}
When $\rho_1=\frac{1}{2}$, we have $\min\{\rho_1,1-\rho_1\}=\frac{1}{2}$. The optimal convergence rate is $O(1/\sqrt{T})$. The proof is completed.
\end{proof}

\begin{rmk}
Our previous works propose a distributed  zeroth-order mirror descent algorithm for constrained optimization over time-varying network in \cite{Yu2,Yuan9,Yuan10} and distributed stochastic mirror descent method for strongly convex objective functions in \cite{Yuan5}. A novel online distributed mirror descent method has been proposed for composite objective functions in \cite{Yuan2}.
In this paper, we largely improve  our previous works \cite{Yuan2, Yu2, Yuan5} on distributed mirror descent methods in several aspects. In the literature of distributed mirror descent, this paper is the first work to propose the adaptive quantization method to address limited communication channel  and simultaneously take  delayed subgradient information into consideration in distributed mirror descent type algorithms. Moreover, the optimal convergence rate $O(1/\sqrt{T})$ is derived under appropriate conditions.
\end{rmk}

\section{Simulations}

In this section, we consider the following distributed estimation problem \cite{Nedic2}:
\begin{eqnarray*}
\min_{x\in\mathcal{X}}\frac{1}{N}\sum^N_{j=1}a_j||x-\textbf{b}_j||^2,
\end{eqnarray*}
over a sequence of time-varying  sensor network, where $a_j\in \mathbb{R}$ and $\textbf{b}_j\in \mathbb{R}^n$. The size of the network is $N=30$ and the dimension of the variable $x$ is $n=10$. The sequence of time-varying network satisfies $B$ connectivity.


 The domain of this system $\mathcal{X}$ is bounded and closed with
\begin{eqnarray*}
\mathcal{X}=\{x\in \mathbb{R}^n|-100\leq [x]_j\leq 100,\ j=1,2,\cdots,n\}.
\end{eqnarray*}
We will use mirror descent algorithm under delayed subgradient information $\tau$, where $\tau=0$ means that subgradient information is timely. We choose the distance generating function $\phi(x)=\frac{1}{2}||x||^2$ and the corresponding Bregman divergence $V_{\phi}(x,y)=\frac{1}{2}||x-y||^2$. Note that $\phi(x)$ satisfies Lipschitz continuity and $V_{\phi}(x,y)$ satisfies separate convexity. Our assumption is satisfied.  In order to show $f(\hat{x}_l(T))$ converges to $f(x^*)$ intuitively, we use the relative error $e(T)$ as
\begin{eqnarray*}
e(T)=\left | \frac{f(\hat{x}_l(T))-f(x^*)}{f(x^*)}\right|
\end{eqnarray*}
\noindent to show the relative error between $f(\hat{x}_l(T))$ and $f(x^*)$.

\subsection{Time Delay's Effect}
\begin{figure}[H]
\centering
\subfigure[without Quantizer]
{
	\begin{minipage}{5cm}
	\centering
	\includegraphics[scale=0.4]{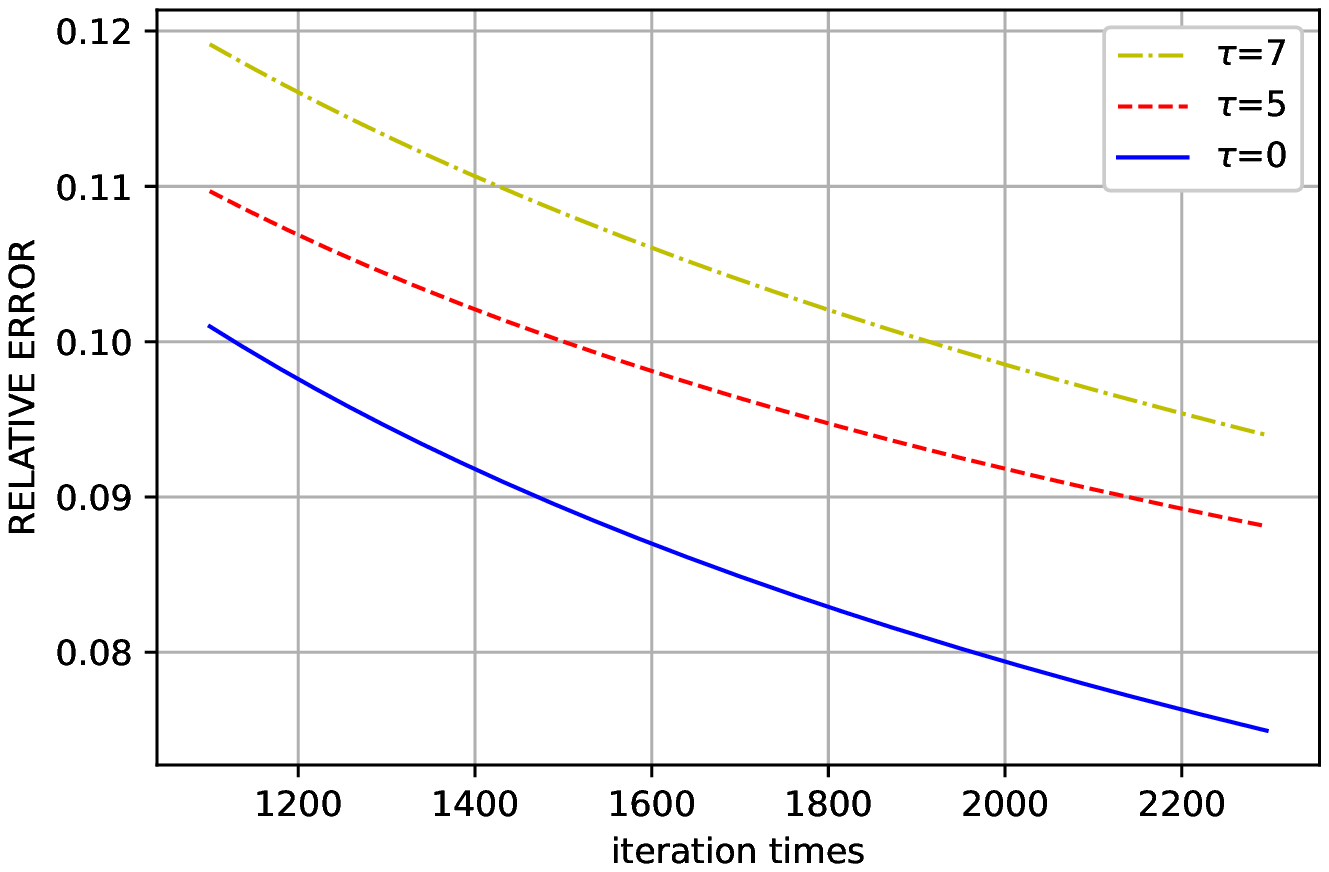}
	\end{minipage}
}
\hspace{.3in}
\subfigure[with Quantizer]
{
	\begin{minipage}{5cm}
	\centering
	\includegraphics[scale=0.4]{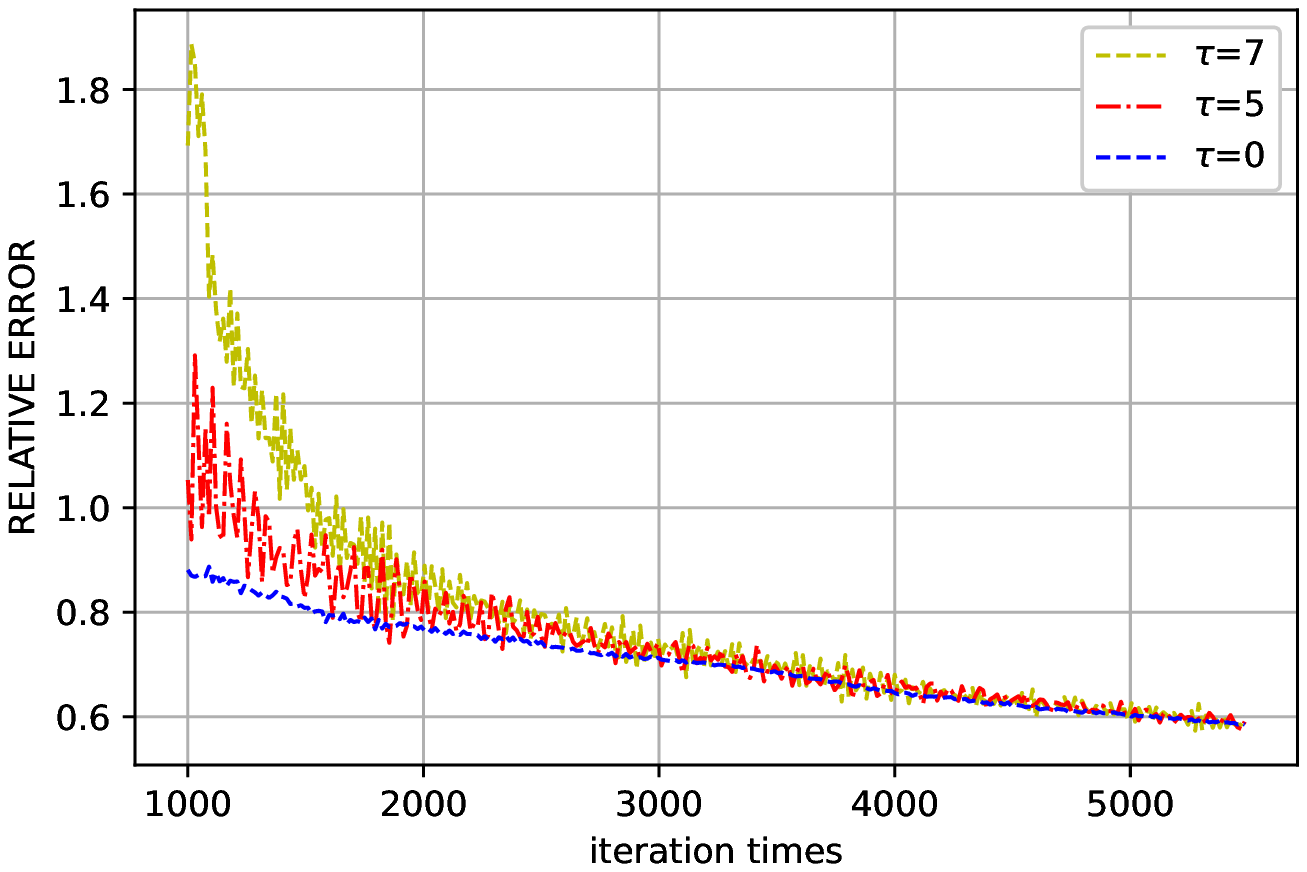}
	\end{minipage}
}
\caption{Time Delay's Effect Comparison}
\label{4}
\end{figure}

In this subsection, we will show the delayed subgradient information's effect to the convergence rate of \textbf{Algorithm 1}. We will choose stepsize $\alpha(t)=\frac{1}{\sqrt{t+1}}$ and quantization parameter $\beta(t)=\frac{1}{\sqrt{t+1}}$ under limited communication channel.  The Fig. \ref{4}(a) shows the convergence under perfect communication channel with time delay $\tau=0,5,7$ and Fig. \ref{4}(b) shows the convergence under limited communication channel with time delay $\tau=0,5,7$. In the Fig. \ref{4}, each line means the relative error between average state of all nodes $\frac{1}{N}\sum^N_{l=1}f(\hat{x}_l(T))$ and optimal value $f(x^*)$.

$\\$ Fig. \ref{4} shows that $\frac{1}{N}\sum^N_{l=1}f(\hat{x}_l(T))$ converges to $f(x^*)$ gradually with time delay $\tau=0$, $\tau=5$ and $\tau=7$ under perfect communication channel and limited communication channel. From Fig. \ref{4}(a), we know that the larger time delay, the slower convergence rate under perfect communication channel. From Fig. \ref{4}(b), the conclusion is similar. The numerical  experiment shows the  larger time delay, the slower convergence rate. In addition, compared with Fig. \ref{4}(a), we can easily find that vibrations of lines in Fig. \ref{4}(b) is much more obvious and the convergence rate is much slower. It is the quantization effect to the convergence.

\subsection{Quantization's Effect}

\begin{figure}[H]
\centering
\subfigure[$\tau=0$]
{
	\begin{minipage}{5cm}
	\centering
	\includegraphics[scale=0.4]{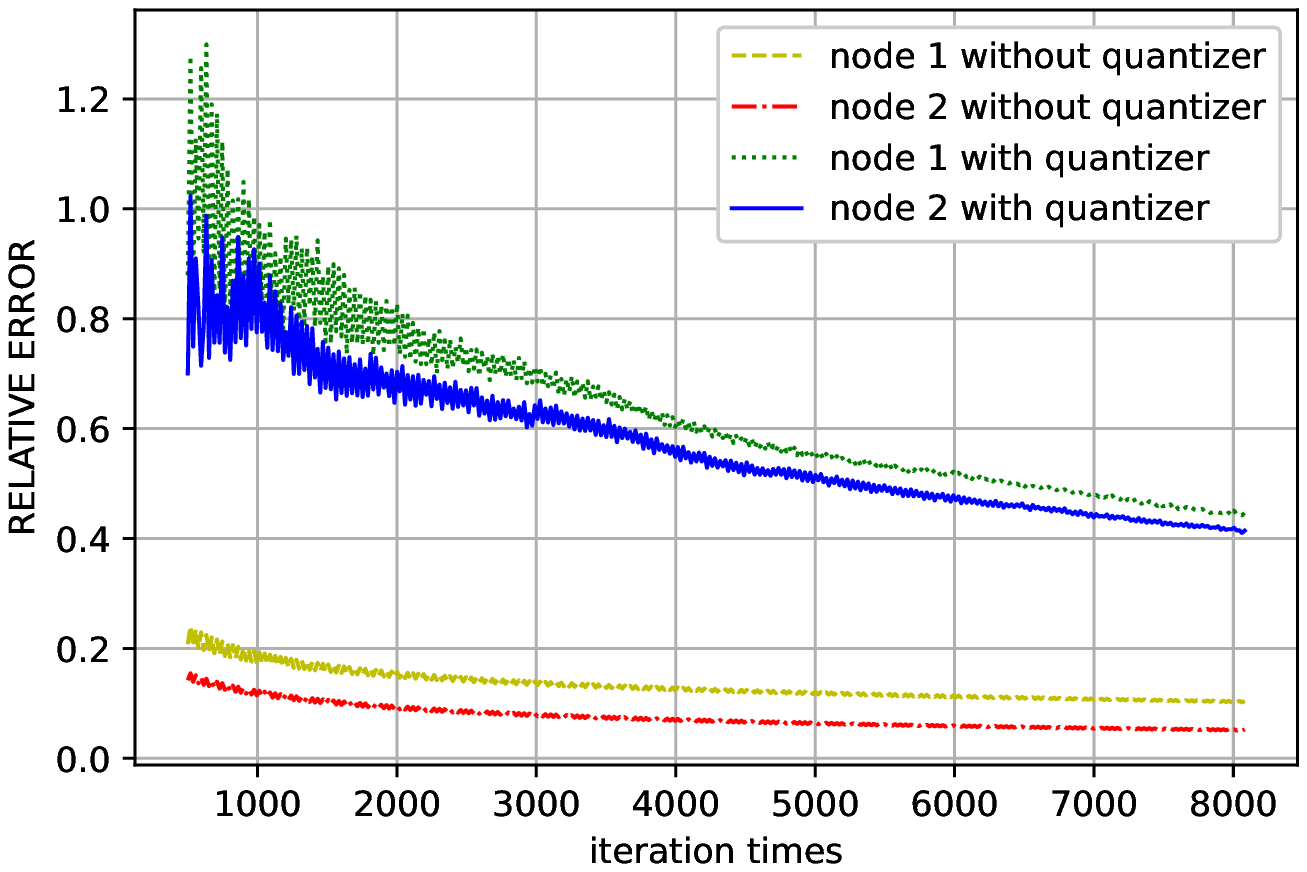}
	\end{minipage}
}
\hspace{.3in}
\subfigure[ $\tau=5$]
{
	\begin{minipage}{5cm}
	\centering
	\includegraphics[scale=0.4]{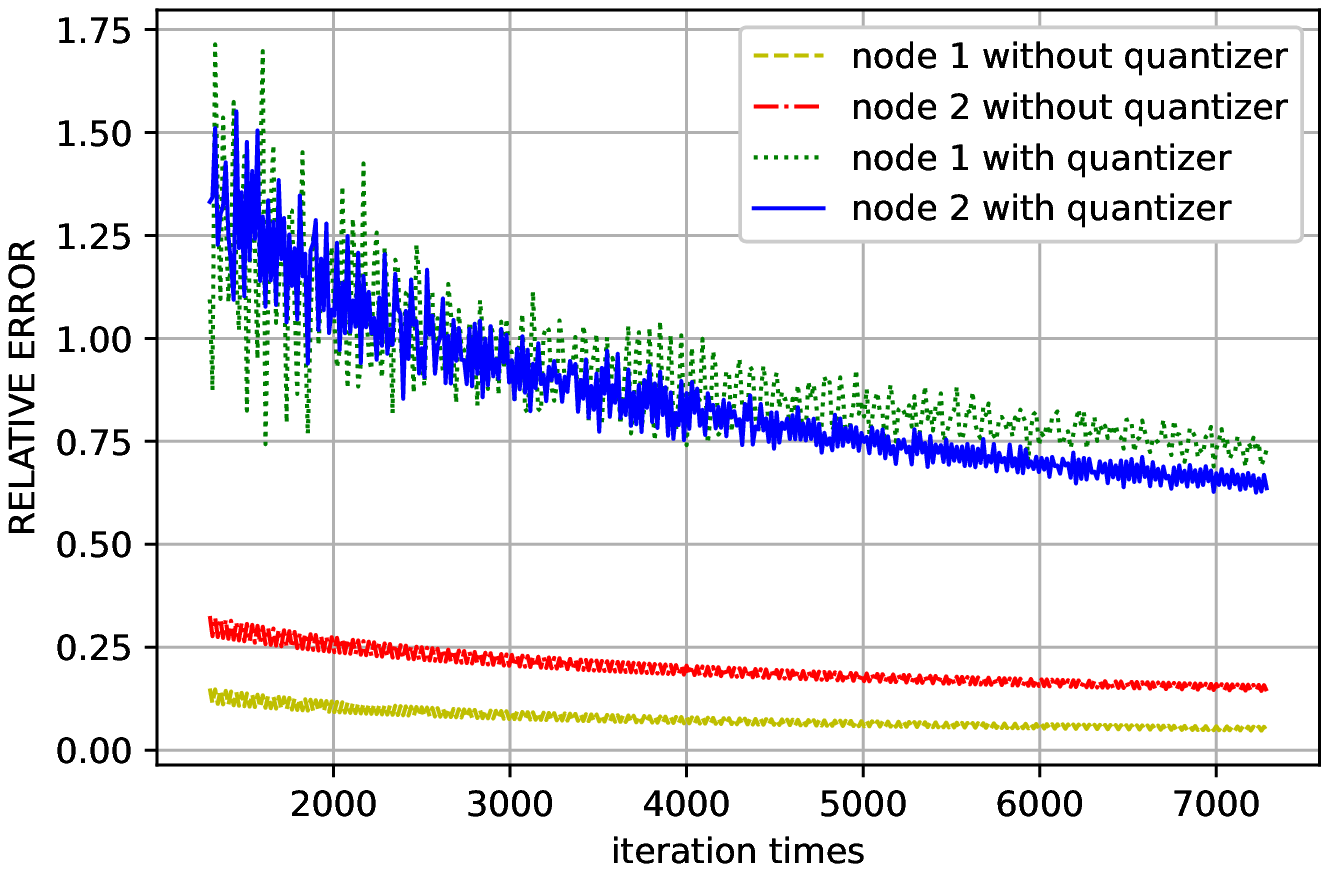}
	\end{minipage}
}
\caption{Quantizer's Effect Comparison}
\label{3}
\end{figure}

In this subsection, we will show the quantization's effect to the convergence of \textbf{Algorithm 1}. We consider  \textbf{Algorithm 1} with stepsize $\alpha(t)=\frac{1}{\sqrt{t+1}}$ under perfect communication channel and limited communication channel ($\beta(t)=\frac{1}{\sqrt{t+1}}$) with different time delay. In the Fig. \ref{3}, we  choose node 1 and node 2 arbitrarily among the network for illustration. The line means the relative error between arbitrarily selected node $f(\hat{x}_l(T))$ and optimal value $f(x^*)$.

$\\$ In Fig. \ref{3}(a),  blue, green lines  show the convergence  under limited communication channel   and red, yellow lines  show the convergence  under perfect communication channel with timely subgradient information ($\tau=0$). In Fig. \ref{3}(b),  blue, green lines  show the convergence   under limited communication channel  and red, yellow lines  show the convergence  under perfect communication channel with delayed subgradient information ($\tau=5$).

From Fig. \ref{3}, we can find that each line converges to $0$ gradually. However,  we can also find that the convergence rate under quantizer is much slower than that under perfect communication channel. Furthermore, we can find that blue and green lines' vibrations are more obvious than those in yellow and red lines' vibrations, which is the quantizer's effect to the convergence of \textbf{Algorithm 1}. The Fig. \ref{3}  reflects quantization will slow down the  convergence rate.

\subsection{Stepsize's Effect}

\begin{figure}[H]
\centering
\includegraphics[width=3in]{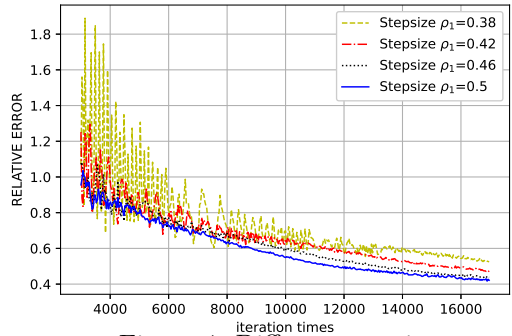}
\caption{Different stepsize}
\label{6}
\end{figure}

In this subsection, we will show the different  stepsize's effect to the convergence rate. We will choose quantization parameter $\beta(t)=\frac{1}{\sqrt{t+1}}$ and different stepsize $\alpha(t)=\frac{1}{(t+1)^{\rho_1}}$, where $\rho_1=0.38,0.42,0.46,0.5$. The following picture shows the convergence of \textbf{Algorithm 1} with different  stepsize. In the Fig. \ref{6}, each line means the relative error between average state of all nodes $\frac{1}{N}\sum^N_{l=1}f(\hat{x}_l(T))$ and optimal value $f(x^*)$.

$\\$ From the above Fig. \ref{6}, we can find that each line converges to $0$ gradually. However, the convergence rate is obviously different. From Corollary \ref{cor22222222}, we know that the correspond convergence rate is $O(1/T^{\rho_1})$ for $\rho_1=0.38,0.42,0.46,0.5$. The numerical  experiment verifies the theoretical result of Corollary \ref{cor22222222}.

\subsection{Quantization Parameter's Effect}

\begin{figure}[H]
\centering
\includegraphics[width=3in]{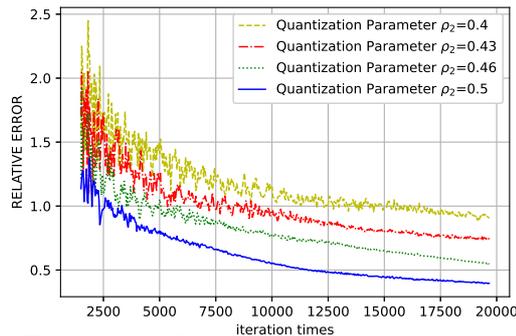}
\caption{Different Quantization Parameter}
\label{5}
\end{figure}

In this subsection, we will show the different quantization parameters' effect to the convergence rate. We will choose stepsize $\alpha(t)=\frac{1}{\sqrt{t+1}}$ and different quantization parameter $\beta(t)=\frac{1}{(t+1)^{\rho_2}}$, where $\rho_2=0.4,0.43,0.46,0.5$. The following picture shows the convergence of \textbf{Algorithm 1} with different quantization parameters. In the Fig. \ref{5}, each line means the relative error between average state of all nodes $\frac{1}{N}\sum^N_{l=1}f(\hat{x}_l(T))$ and optimal value $f(x^*)$.

$\\$ From the above Fig. \ref{5}, we can find that each line converges to $0$ gradually. However the convergence rate is obviously different. We can conclude from Fig. \ref{5} that the larger $\rho_2$, the faster it converges. From  Corollary \ref{cor22222222}, we know that the correspond convergence rate is $O(1/T^{\rho_2})$  for $\rho_2=0.4,0.43,0.46,0.5$. The numerical  experiment verifies the theoretical result of Corollary \ref{cor22222222}.

\section{Conclusion}
This paper has studied the mirror descent algorithm with delayed subgradient information under adaptive quantization.  We design adaptive quantization method to solve the limited communication channel and analyze the convergence of \textbf{Algorithm 1} with different stepsizes and quantization parameters. The optimal convergence rate $O(\frac{1}{\sqrt{T}})$ can be obtained under appropriate conditions. This paper has improved many previous works (e.g. \cite{Li,Yuan2, Yu2, Yuan5,Doan}) in many aspects. Some numerical examples have been presented to demonstrate the effectiveness of the algorithm and verify the theoretical results. For the future work, we can consider how to use adaptive method in other settings, such as dual averaging algorithm \cite{Duchi1}, online composite optimization \cite{Yuan2,Lee}.



\section{Appendix}

\subsection{Proof of Lemma \ref{lemma-6}}

Before proofing Lemma \ref{lemma-6}, we need the following Lemma \ref{lemma-3} and  Lemma \ref{lemma-5}.

\begin{lem}\label{lemma-3}
$\{\tilde{y}_i(t)\}$ and $\{g_i(t-\tau)\}$ are sequences generated by \textbf{Algorithm 1}. Then we have
\begin{eqnarray*}
\langle g_i(t-\tau),\tilde{y}_i(t)-x^*\rangle\leq \frac{1}{\alpha(t)}V_{\phi}(x^*,\tilde{y}_i(t))-\frac{1}{\alpha(t)}V_{\phi}(x^*,x_i(t+1))+\frac{G^2\alpha(t)}{2\sigma_{\phi}}.
\end{eqnarray*}
\end{lem}


\begin{proof}
For the first order optimality condition, for $\forall x\in \mathcal{X}$ we have
\begin{eqnarray*}
\langle \alpha(t)g_i(t-\tau)+\nabla \phi(x_i(t+1))-\nabla \phi(\tilde{y}_i(t)),x-x_i(t+1)\rangle \geq 0
\end{eqnarray*}
\noindent Thus we select $x=x^*$ and we have
\begin{eqnarray}\label{abcdefg}
\langle\nabla \phi(x_i(t+1))-\nabla \phi(\tilde{y}_i(t)),x^*-x_i(t+1)\rangle\geq \langle\alpha(t)g_i(t-\tau),x_i(t+1)-x^*\rangle.
\end{eqnarray}

\noindent Rearrange terms in (\ref{abcdefg}) and we have
\begin{eqnarray}
&&\langle\alpha(t)g_i(t-\tau),x_i(t+1)-x^*\rangle\nonumber\\
&\leq& \langle\nabla \phi(\tilde{y}_i(t))-\nabla \phi(x_i(t+1)),x_i(t+1)-x^*\rangle\label{ddffddff}\\
&=&V_{\phi}(x^*,\tilde{y}_i(t))-V_{\phi}(x^*,x_i(t+1))-V_{\phi}(x_i(t+1),\tilde{y}_i(t))\nonumber\\
&\leq&V_{\phi}(x^*,\tilde{y}_i(t))-V_{\phi}(x^*,x_i(t+1))-\frac{\sigma_{\phi}}{2}||x_i(t+1)-\tilde{y}_i(t)||^2,\label{lemma31}
\end{eqnarray}
where the  inequality (\ref{ddffddff}) follows from Bregman divergence's definition and  the  inequality (\ref{lemma31}) follows from the strong convexity of $\phi$. From Cauchy inequality, we have
\begin{eqnarray}\label{lemma32}
&&\langle\alpha(t)g_i(t-\tau),x_i(t+1)-x^*\rangle\nonumber\\
&=&\langle\alpha(t)g_i(t-\tau),x_i(t+1)-\tilde{y}_i(t)\rangle+\langle\alpha(t)g_i(t-\tau),\tilde{y}_i(t)-x^*\rangle\nonumber\\
&\geq&-\frac{\alpha^2(t)}{2\sigma_{\phi}}||g_i(t-\tau)||^2-\frac{\sigma_{\phi}}{2}||x_i(t+1)-\tilde{y}_i(t)||^2+\langle\alpha(t)g_i(t-\tau),\tilde{y}_i(t)-x^*\rangle.
\end{eqnarray}
Combining inequality (\ref{lemma31}) and (\ref{lemma32}), we have
\begin{eqnarray*}
&&\langle g_i(t-\tau),\tilde{y}_i(t)-x^*\rangle\\
&\leq& \frac{\alpha(t)}{2\sigma_{\phi}}||g_i(t-\tau)||^2+\frac{1}{\alpha(t)}\{V_{\phi}(x^*,\tilde{y}_i(t))-V_{\phi}(x^*,x_i(t+1))\}\\
&\leq& \frac{G^2\alpha(t)}{2\sigma_{\phi}}+\frac{1}{\alpha(t)}\{V_{\phi}(x^*,\tilde{y}_i(t))-V_{\phi}(x^*,x_i(t+1))\}.
\end{eqnarray*}

\noindent The proof of Lemma \ref{lemma-2} is completed.
\end{proof}

\begin{lem}\label{lemma-5}
$\{x_j(t)\}$ is generated by  \textbf{Algorithm 1}. $e_j(t)$ and $p_i(t)$ are defined as equations (\ref{quantizationerror}) and (\ref{projectionerror}). Then we have
\begin{eqnarray*}
&&V_{\phi}(x^*,x_j(t)+e_j(t)+p_i(t))\leq V_{\phi}(x^*,x_j(t))+2(4N^2+1)L_{\phi}E^2(t)+(2N+1)L_{\phi}\sqrt{\frac{2D_{\phi}}{\sigma_{\phi}}}E(t).
\end{eqnarray*}
\end{lem}

\begin{proof}
\noindent Follow from mean value formula, there exists a $0\leq \xi\leq 1$ such that
\begin{eqnarray*}
\phi(x_j(t)+e_j(t)+p_i(t))=\phi(x_j(t))+\langle \nabla\phi(x_j(t)+\xi e_j(t)+\xi p_i(t)),e_j(t)+p_i(t)\rangle.
\end{eqnarray*}

\noindent Then we have
\begin{eqnarray}\label{lemma51}
&&V_{\phi}(x^*,x_j(t)+e_j(t)+p_i(t))\nonumber\\
&=&\phi(x^*)-\phi(x_j(t)+e_j(t)+p_i(t))-\langle \nabla(x_j(t)+e_j(t)+p_i(t),x^*-x_j(t)-e_j(t)-p_i(t))\rangle\nonumber\\
&=&\phi(x^*)-\phi(x_j(t))-\langle \nabla(x_j(t)+e_j(t)+p_i(t),x^*-x_j(t)-e_j(t)-p_i(t))\rangle\nonumber\\
&&-\langle \nabla\phi(x_j(t)+\xi e_j(t)+\xi p_i(t)),e_j(t)+p_i(t)\rangle\nonumber\\
&=&\phi(x^*)-\phi(x_j(t))-\langle \nabla\phi(x_j(t)), x^*-x_j(t)\rangle\nonumber\\
&&-\langle \nabla\phi(x_j(t)+e_j(t)+p_i(t))-\nabla\phi(x_j(t)), x^*-x_j(t)\rangle\nonumber\\
&&-\langle \nabla\phi(x_j(t)+\xi e_j(t)+\xi p_i(t))- \nabla\phi(x_j(t)+e_j(t)+p_i(t)),e_j(t)+p_i(t) \rangle\nonumber\\
&\leq&V_{\phi}(x^*,x_j(t))+(1-\xi)L_{\phi}||e_j(t)+p_i(t)||^2+L_{\phi}\sqrt{\frac{2D_{\phi}}{\sigma_{\phi}}}||e_j(t)+p_i(t)||,
\end{eqnarray}
where the inequality follows from Lipschitz continuous of $\phi$. We note that
\begin{eqnarray*}
||e_j(t)+p_i(t)||\leq||p_i(t)||+||e_j(t)||\leq (2N+1)E(t)\quad\quad\quad&\\
||e_j(t)+p_i(t)||^2\leq2||e_j(t)||^2+2||p_i(t)||^2\leq 2E^2(t)+8N^2E^2(t).&
\end{eqnarray*}
Substitute the above two inequalities into inequality (\ref{lemma51}), we have
\begin{eqnarray*}
V_{\phi}(x^*,x_j(t)+e_j(t)+p_i(t))&\leq& V_{\phi}(x^*,x_j(t))+2(1+4N^2)L_{\phi}E^2(t)\\
&&+(2N+1)L_{\phi}\sqrt{\frac{2D_{\phi}}{\sigma_{\phi}}}E(t).
\end{eqnarray*}
 The proof of Lemma \ref{lemma-5} is completed.
\end{proof}
Then, the proof of Lemma \ref{lemma-6} is shown as follows.
\begin{proof}
Follow from Lemma \ref{lemma-3}, we have
\begin{eqnarray}\label{lemma661}
&&\sum^{T+\tau}_{t=1+\tau}\sum^N_{i=1}\langle g_i(t-\tau),\tilde{y}_i(t)-x^*\rangle\nonumber\\
&\leq&\sum^{T+\tau}_{t=1+\tau}\frac{1}{\alpha(t)}\sum^N_{i=1}\{V_{\phi}(x^*,\tilde{y}_i(t))-V_{\phi}(x^*,x_i(t+1))\}+\frac{NG^2}{2\sigma_{\phi}}\sum^{T+\tau}_{t=1+\tau}\alpha(t).
\end{eqnarray}
For the
\begin{eqnarray*}
\sum^T_{t=1}\frac{1}{\alpha(t)}\sum^N_{i=1}\{V_{\phi}(x^*,\tilde{y}_i(t))-V_{\phi}(x^*,x_i(t+1))\},
\end{eqnarray*}
we have
\begin{eqnarray}
&&\sum^{T+\tau}_{t=1+\tau}\frac{1}{\alpha(t)}\sum^N_{i=1}\Big\{V_{\phi}(x^*,\tilde{y}_i(t))-V_{\phi}(x^*,x_i(t+1))\Big\}\nonumber\\
&=&\sum^{T+\tau}_{t=1+\tau}\frac{1}{\alpha(t)}\Big\{\sum^N_{i=1}V_{\phi}(x^*,\tilde{y}_i(t))-\sum^N_{i=1}V_{\phi}(x^*,x_i(t+1))\Big\}\nonumber\\
&\leq&\sum^{T+\tau}_{t=1+\tau}\frac{1}{\alpha(t)}\Big\{\sum^N_{i=1}\sum^N_{j=1}[P(t)]_{ij}V_{\phi}(x^*,x_j(t)+e_j(t)+p_i(t))-\nonumber\\
&&\sum^N_{i=1}\sum^N_{j=1}[P(t)]_{ij}V_{\phi}(x^*,x_i(t+1))\Big\}\label{666662}\\
&=&\sum^{T+\tau}_{t=1+\tau}\frac{1}{\alpha(t)}\Big\{\sum^N_{i=1}\sum^N_{j=1}[P(t)]_{ij}\{V_{\phi}(x^*,x_j(t)+e_j(t)+p_i(t))-V_{\phi}(x^*,x_i(t+1))\}\Big\}\nonumber\\
&\leq&\sum^{T+\tau}_{t=1+\tau}\frac{1}{\alpha(t)}\Big\{\sum^N_{i=1}\sum^N_{j=1}[P(t)]_{ij}\{V_{\phi}(x^*,x_j(t))-V_{\phi}(x^*,x_i(t+1))\}\Big\}\nonumber\\
&&+\sum^{T+\tau}_{t=1+\tau}\frac{1}{\alpha(t)}\Big\{\sum^N_{i=1}\sum^N_{j=1}[P(t)]_{ij}\{2(1+4N^2)L_{\phi}E^2(t)+(2N+1)L_{\phi}\sqrt{\frac{2D_{\phi}}{\sigma_{\phi}}}E(t)\} \Big\}\label{666663}\\
&=&\sum^{T+\tau}_{t=1+\tau}\frac{1}{\alpha(t)}\Big\{\sum^N_{j=1}V_{\phi}(x^*,x_j(t))-\sum^N_{i=1}V_{\phi}(x^*,x_i(t+1)) \Big\}+L_{\phi}\sqrt{\frac{2D_{\phi}}{\sigma_{\phi}}}(2N^2+N)\sum^{T+\tau}_{t=1+\tau}\frac{E(t)}{\alpha(t)}\nonumber\\
&&+2L_{\phi}(4N^3+N)\sum^{T+\tau}_{t=1+\tau}\frac{E^2(t)}{\alpha(t)},\label{lemma662}
\end{eqnarray}
where the first inequality (\ref{666662}) follows from separate convexity of $V_{\phi}(x,y)$ and second inequality (\ref{666663}) follows from Lemma \ref{lemma-5}.
For the $\sum^{T+\tau}_{t=1+\tau}\frac{1}{\alpha(t)}\Big\{\sum^N_{j=1}V_{\phi}(x^*,x_j(t))-\sum^N_{i=1}V_{\phi}(x^*,x_i(t+1)) \Big\}$, we have
\begin{eqnarray}
&&\sum^{T+\tau}_{t=1+\tau}\frac{1}{\alpha(t)}\Big\{\sum^N_{j=1}V_{\phi}(x^*,x_j(t))-\sum^N_{i=1}V_{\phi}(x^*,x_i(t+1)) \Big\}\nonumber\\
&=&\sum^N_{i=1}\Big\{\sum^{T+\tau}_{t=2+\tau}V_{\phi}(x^*,x_j(t))(\frac{1}{\alpha(t)}-\frac{1}{\alpha(t-1)})+\frac{1}{\alpha(1+\tau)}V_{\phi}(x^*,x_j(1+\tau))\nonumber\\
&&-\frac{1}{\alpha(T+\tau)}V_{\phi}(x^*,x_j(T+1+\tau))  \Big\}\nonumber\\
&\leq&\sum^N_{j=1}D_{\phi}\{\frac{1}{\alpha(1+\tau)}+\sum^{T+\tau}_{t=2+\tau}(\frac{1}{\alpha(t)}-\frac{1}{\alpha(t-1)})\}\label{lemma66553}\\
&=&\frac{ND_{\phi}}{\alpha(T+\tau)},\label{lemma663}
\end{eqnarray}
where the first inequality (\ref{lemma66553}) follows from the upper bound of $V_{\phi}(x,y)$.  Substitute inequality (\ref{lemma663}) into (\ref{lemma662}) and then  substitute inequality (\ref{lemma662}) into (\ref{lemma661}), we have
\begin{eqnarray*}
&&\sum^{T+\tau}_{t=1+\tau}\sum^N_{i=1}\langle g_i(t-\tau),\tilde{y}_i(t)-x^*\rangle\\
&\leq&(2N^2+N)L_{\phi}\sqrt{\frac{2D_{\phi}}{\sigma_{\phi}}}\sum^{T+\tau}_{t=1+\tau}\frac{E(t)}{\alpha(t)}+\frac{N G^2}{2\sigma_{\phi}}\sum^{T+\tau}_{t=1+\tau}\alpha(t)+\frac{ND_{\phi}}{\alpha(T+\tau)}+(8N^3+2N)L_{\phi}\sum^{T+\tau}_{t=1+\tau}\frac{E^2(t)}{\alpha(t)}.
\end{eqnarray*}
\noindent Therefore, we have
\begin{eqnarray*}
&&\frac{1}{NT}\sum^{T+\tau}_{t=1+\tau}\sum^N_{i=1}\langle g_i(t-\tau),\tilde{y}_i(t)-x^*\rangle\\
&\leq&(2N+1)L_{\phi}\sqrt{\frac{2D_{\phi}}{\sigma_{\phi}}}\frac{1}{T}\sum^{T+\tau}_{t=1+\tau}\frac{E(t)}{\alpha(t)}+\frac{G^2}{2\sigma_{\phi}}\frac{1}{T}\sum^{T+\tau}_{t=1+\tau}\alpha(t)+\frac{D_{\phi}}{T\alpha(T+\tau)}+L_{\phi}(8N^2+2)\frac{1}{T}\sum^{T+\tau}_{t=1+\tau}\frac{E^2(t)}{\alpha(t)}\\
&=&\mathcal{B}_1
\end{eqnarray*}
The proof is completed.
\end{proof}
\subsection{Proof of Lemma \ref{lemma-7}}
\begin{proof}
From Cauchy inequality and the upper bound of subgradient, we have
\begin{eqnarray}\label{lemma771}
\langle g_i(t-\tau),\tilde{y}_i(t-\tau)-\tilde{y}_i(t)\rangle\leq||g_i(t-\tau)||||\tilde{y}_i(t-\tau)-\tilde{y}_i(t)||\nonumber\\
&\leq&G||\tilde{y}_i(t-\tau)-\tilde{y}_i(t)||.
\end{eqnarray}
After summing up inequality (\ref{lemma771}) from $t =\tau+ 1$ to $T + \tau$ and dividing both side by $NT$, we have
\begin{eqnarray}\label{lemma772}
\frac{1}{NT}\sum^{T+\tau}_{t=1+\tau}\sum^N_{i=1}\langle g_i(t-\tau),\tilde{y}_i(t-\tau)-\tilde{y}_i(t)\rangle\leq \frac{G}{NT}\sum^{T+\tau}_{t=1+\tau}\sum^N_{i=1}||\tilde{y}_i(t-\tau)-\tilde{y}_i(t)||.
\end{eqnarray}
For  $||\tilde{y}_i(t-\tau)-\tilde{y}_i(t)||$, from triangle inequality, we have
\begin{eqnarray}\label{773}
||\tilde{y}_i(t-\tau)-\tilde{y}_i(t)||\leq\sum^{t-1}_{k=t-\tau}||\tilde{y}_i(k)-\tilde{y}_i(k+1)||.
\end{eqnarray}
For $||\tilde{y}_i(k)-\tilde{y}_i(k+1)||$, we have
\begin{eqnarray}
&&||\tilde{y}_i(k)-\tilde{y}_i(k+1)||\nonumber\\
&\leq&||\tilde{y}_i(k)-x_i(k+1)||+||x_i(k+1)-\tilde{y}_i(k+1)||\label{77664}\\
&\leq&||\sum^N_{j=1}[P(k+1)]_{ij}(x_j(k+1))-x_i(k+1)||+||\varepsilon_i(k)||\nonumber\\
&&+||\sum^N_{j=1}[P(k+1)]_{ij}(e_j(k+1)+p_i(k+1))||\label{77665}\\
&\leq&\frac{G\alpha(k)}{\sigma_{\phi}}+3NE(k+1)+\sum^N_{j=1}||x_j(k+1)-x_i(k+1)||,\label{774}
\end{eqnarray}
where inequality (\ref{77664}) and  inequality  (\ref{77665})  follow from triangle inequality. The sequence $\{\alpha(t)\}$ and $\{E(t)\}$ are non-increasing and substitute inequality (\ref{774}) into (\ref{773}), so we have
\begin{eqnarray}\label{775}
&&||\tilde{y}_i(t-\tau)-\tilde{y}_i(t)||\nonumber\\
&\leq&\sum^{t-1}_{k=t-\tau}||\tilde{y}_i(k)-\tilde{y}_i(k+1)||\nonumber\\
&\leq&\sum^{t-1}_{k=t-\tau}\big\{\frac{G\alpha(k)}{\sigma_{\phi}}+3NE(k+1)+\sum^N_{j=1}||x_j(k+1)-x_i(k+1)||\big\}\nonumber\\
&\leq&\frac{G\tau}{\sigma_{\phi}}\alpha(t-\tau)+3N\tau E(t-\tau)+\sum^{t-1}_{k=t-\tau}\sum^N_{j=1}||x_j(k+1)-x_i(k+1)||.
\end{eqnarray}
Substitute inequality (\ref{775}) into (\ref{lemma772}) and we have
\begin{eqnarray}\label{7777}
&&\frac{G}{NT}\sum^{T+\tau}_{t=1+\tau}\sum^N_{i=1}||\tilde{y}_i(t-\tau)-\tilde{y}_i(t)||\nonumber\\
&\leq& \frac{G}{NT}\sum^{T+\tau}_{t=1+\tau}\sum^N_{i=1}\big\{\frac{G\tau}{\sigma_{\phi}}\alpha(t-\tau)+3N\tau E(t-\tau)+\sum^{t-1}_{k=t-\tau}\sum^N_{j=1}||x_j(k+1)-x_i(k+1)||\big\}\nonumber\\
&=&\frac{G}{NT}\sum^{T+\tau}_{t=1+\tau}\sum^N_{i=1}\sum^{t-1}_{k=t-\tau}\sum^N_{j=1}||x_j(k+1)-x_i(k+1)||+\frac{G^2\tau}{\sigma_{\phi}T}\sum^{T+\tau}_{t=1+\tau}\alpha(t-\tau)\nonumber\\
&&+\frac{3N\tau G}{T}\sum^{T+\tau}_{t=1+\tau} E(t-\tau).
\end{eqnarray}
For $\frac{G}{NT}\sum^{T+\tau}_{t=1+\tau}\sum^N_{i=1}\sum^{t-1}_{k=t-\tau}\sum^N_{j=1}||x_j(k+1)-x_i(k+1)||$, we have
\begin{eqnarray}\label{776}
&&\frac{G}{NT}\sum^{T+\tau}_{t=1+\tau}\sum^N_{i=1}\sum^{t-1}_{k=t-\tau}\sum^N_{j=1}||x_j(k+1)-x_i(k+1)||\nonumber\\
&=&\frac{G}{NT}\sum^{T+\tau}_{t=1+\tau}\sum^N_{i=1}\sum^N_{j=1}\Big\{||x_j(t-\tau+1)-x_i(t-\tau+1)||+||x_j(t-\tau+2)-x_i(t-\tau+2)||\nonumber+\cdots\\
&&+||x_j(t)-x_i(t)||\Big\}\nonumber\\
&=&\frac{G}{NT}\sum^{T}_{t=1}\sum^N_{i=1}\sum^N_{j=1}\Big\{||x_j(t+1)-x_i(t+1)||+||x_j(t+2)-x_i(t+2)||+\cdots\nonumber\\
&&+||x_j(t+\tau)-x_i(t+\tau)||\Big\}\nonumber\\
&\leq&\frac{G\tau}{NT}\sum^{T+\tau}_{t=1}\sum^N_{i=1}\sum^N_{j=1}||x_j(t)-x_i(t)||.
\end{eqnarray}
From Lemma \ref{lemma-2}, we have
\begin{eqnarray*}
\sum^T_{t=1}\sum^N_{i=1}||x_i(t)-x_j(t)||\leq \frac{2N\omega}{1-\gamma}A+2B\sum^T_{t=0}\{\frac{G\alpha(t)}{\sigma_{\phi}}+3NE(t)\}
\end{eqnarray*}
Thus we have
\begin{eqnarray}\label{777}
&&\sum^{T+\tau}_{t=1}\sum^N_{i=1}\sum^N_{j=1}||x_j(t)-x_i(t)||\nonumber\\
&\leq&\frac{2N^2\omega}{1-\gamma}A+2BN\sum^{T+\tau}_{t=0}\{\frac{G\alpha(t)}{\sigma_{\phi}}+3NE(t)\}\nonumber\\
&=&\frac{2N^2\omega}{1-\gamma}A+\frac{2BNG}{\sigma_{\phi}}\sum^{T+\tau}_{t=0}\alpha(t)+6BN^2\sum^{T+\tau}_{t=0}E(t).
\end{eqnarray}
Substitute inequality (\ref{777}) into (\ref{776}) and we have
\begin{eqnarray}\label{888}
&&\frac{G}{NT}\sum^{T+\tau}_{t=1+\tau}\sum^N_{i=1}\sum^{t-1}_{k=t-\tau}\sum^N_{j=1}||x_j(k+1)-x_i(k+1)||\nonumber\\
&\leq& \frac{G\tau}{NT}\big\{\frac{2N^2\omega}{1-\gamma}A+\frac{2BNG}{\sigma_{\phi}}\sum^{T+\tau}_{t=0}\alpha(t)+6BN^2\sum^{T+\tau}_{t=0}E(t)\big\}\nonumber\\
&\leq&\frac{2N G\tau\omega A}{(1-\gamma)T}+\frac{2BG^2\tau}{\sigma_{\phi}T}\sum^{T+\tau}_{t=0}\alpha(t)+\frac{6BN G\tau}{T}\sum^{T+\tau}_{t=0}E(t).
\end{eqnarray}
Substitute inequality (\ref{888}) into (\ref{7777}) and substitute (\ref{7777}) into (\ref{lemma772}), then we have
\begin{eqnarray*}
&&\frac{1}{NT}\sum^{T+\tau}_{t=1+\tau}\sum^N_{i=1}\langle g_i(t-\tau),\tilde{y}_i(t-\tau)-\tilde{y}_i(t)\rangle\\
&\leq&\frac{G}{NT}\sum^{T+\tau}_{t=1+\tau}\sum^N_{i=1}||\tilde{y}_i(t-\tau)-\tilde{y}_i(t)||\\
&\leq&\frac{G^2\tau}{\sigma_{\phi}}\frac{1}{T}\sum^{T+\tau}_{t=1+\tau}\alpha(t-\tau)+\frac{3N\tau G}{T}\sum^{T+\tau}_{t=1+\tau} E(t-\tau)+\frac{2N G\tau\omega A}{(1-\gamma)T}+\frac{2BG^2\tau}{\sigma_{\phi}T}\sum^{T+\tau}_{t=0}\alpha(t)+\frac{6BN G\tau}{T}\sum^{T+\tau}_{t=0}E(t)\\
&=&\frac{G^2\tau}{\sigma_{\phi}}\frac{1}{T}\sum^{T}_{t=1}\alpha(t)+\frac{3N\tau G}{T}\sum^{T}_{t=1} E(t)+\frac{1}{T}\frac{2N G\tau\omega}{1-\gamma}A+\frac{2BG^2\tau}{\sigma_{\phi}}\frac{1}{T}\sum^{T+\tau}_{t=0}\alpha(t)+6BN G\tau\frac{1}{T}\sum^{T+\tau}_{t=0}E(t)\\
&\leq&\frac{(2B+1)G^2\tau}{\sigma_{\phi}}\frac{1}{T}\sum^{T+\tau}_{t=0}\alpha(t)+\frac{(6B+3)N G\tau}{T}\sum^{T+\tau}_{t=0}E(t)+\frac{1}{T}\frac{2N G\tau\omega}{1-\gamma}A\\
&=&\mathcal{B}_2.
\end{eqnarray*}
The proof is completed.
\end{proof}

\subsection{Proof of Lemma \ref{lemma-2}}
\begin{proof}
We define that $A=\sum^N_{j=1}||x_j(0)||$, $B=(2N+\frac{N^2\omega}{1-\gamma})$ and  we will show the bound of $||x_i(t)-\bar{x}(t)||$.
\begin{eqnarray}\label{anotherbound}
||c_j(t)||&=&||\sum^N_{j=1}[P(t)]_{ij}(e_j(t)+p_i(t))||\nonumber\\
&\leq&\sum^N_{j=1}||e_j(t)||+||p_i(t)||\nonumber\\
&\leq&3NE(t).
\end{eqnarray}
 Therefore, we have
\begin{eqnarray*}
&&||x_i(t)-\bar{x}(t)||\\
&\leq&\sum^t_{s=1}\sum^N_{j=1}\big|[P(t-1,s)]_{ij}-\frac{1}{N}\big|||c_j(s-1)+\varepsilon_j(s-1)||+\sum^N_{j=1}\big|[P(t-1,0)]_{ij}-\frac{1}{N}\big|||x_j(0)||\\
&\leq&\omega\gamma^{t-1}A+\sum^{t-1}_{s=1}\omega\gamma^{t-s-1}\sum^N_{j=1}||c_j(s-1)+\varepsilon_j(s-1)||+\frac{1}{N}\sum^N_{j=1}||c_j(t-1)+\varepsilon_j(t-1)||\\
&&+||c_i(t-1)+\varepsilon_i(t-1)||.
\end{eqnarray*}
From the triangle inequality, Lemma \ref{Bregman} and  (\ref{anotherbound}), we have
\begin{eqnarray*}
||\varepsilon_j(s)+c_j(s)||\leq \frac{G\alpha(s)}{\sigma_{\phi}}+3NE(s).
\end{eqnarray*}
\noindent Thus we have
\begin{eqnarray*}
&&||x_i(t)-\bar{x}(t)||\\
&\leq& \omega\gamma^{t-1}A+2(\frac{G\alpha(t-1)}{\sigma_{\phi}}+3NE(t-1))+N\sum^{t-1}_{s=1}\omega\gamma^{t-s-1}(\frac{G\alpha(s-1)}{\sigma_{\phi}}+3NE(s-1)).
\end{eqnarray*}
Then we have
\begin{eqnarray*}
&&||x_i(t)-x_j(t)||\\
&\leq& ||x_j(t)-\bar{x}(t)||+||x_i(t)-\bar{x}(t)||\\
&\leq& 2\omega\gamma^{t-1}A+4(\frac{G\alpha(t-1)}{\sigma_{\phi}}+3NE(t-1))+2N\sum^{t-1}_{s=1}\omega\gamma^{t-s-1}(\frac{G\alpha(s-1)}{\sigma_{\phi}}+3NE(s-1)).
\end{eqnarray*}
At last, we will show the bound of $\sum^T_{t=1}\sum^N_{i=1}||x_i(t)-\bar{x}(t)||$.
\begin{eqnarray*}
\sum^T_{t=1}\sum^N_{i=1}||x_i(t)-\bar{x}(t)||\leq \frac{N\omega}{1-\gamma}A+B\sum^T_{t=0}\{\frac{G\alpha(t)}{\sigma_{\phi}}+3NE(t)\}
\end{eqnarray*}

\begin{eqnarray*}
\sum^T_{t=1}\sum^N_{i=1}||x_i(t)-x_j(t)||\leq \frac{2N\omega}{1-\gamma}A+2B\sum^T_{t=0}\{\frac{G\alpha(t)}{\sigma_{\phi}}+3NE(t)\}.
\end{eqnarray*}
The proof is completed.
\end{proof}
\subsection{Proof of Corollary \ref{cor6666} }
\begin{proof}
Obviously, we know that
\begin{eqnarray*}
\lim_{T\rightarrow \infty}\frac{1}{T}\sum^{T}_{t=1}\alpha(t)=0&\Leftrightarrow& \lim_{T\rightarrow \infty}\frac{1}{T}\sum^{T+\tau}_{t=1}\alpha(t)=0,\\
\lim_{T\rightarrow \infty}\frac{1}{T}\sum^{T}_{t=1}\beta(t)=0&\Leftrightarrow& \lim_{T\rightarrow \infty}\frac{1}{T}\sum^{T+\tau}_{t=1}\beta(t)=0.
\end{eqnarray*}
When $\lim_{T\rightarrow \infty}\frac{1}{T\alpha(T+\tau)}=0$, we have
\begin{eqnarray*}
\frac{1}{T\alpha(T)}\leq \frac{1}{T\alpha(T+\tau)}.
\end{eqnarray*}
Therefore, we have $\lim_{T\rightarrow \infty}\frac{1}{T\alpha(T)}=0$.

\noindent When $\lim_{T\rightarrow \infty}\frac{1}{T\alpha(T)}=0$, we have
\begin{eqnarray*}
&&\lim_{T\rightarrow \infty}\frac{1}{T\alpha(T+\tau)}\\
&=&\lim_{T\rightarrow \infty}\frac{T+\tau}{T(T+\tau)\alpha(T+\tau)}\\
&=&\lim_{T\rightarrow \infty}\frac{1}{(T+\tau)\alpha(T+\tau)}\lim_{T\rightarrow \infty}\frac{T+\tau}{T}\\
&=&\lim_{T+\tau\rightarrow \infty}\frac{1}{(T+\tau)\alpha(T+\tau)}\lim_{T\rightarrow \infty}\frac{T+\tau}{T}\\
&=&\lim_{T\rightarrow \infty}\frac{1}{T\alpha(T)}\\
&=&0.
\end{eqnarray*}
Therefore, we can conclude that  (\ref{AAAA2}), (\ref{BBBB2}), (\ref{CCCCdsfsad2}) are equivalent to (\ref{DDDD2}), (\ref{EEEE2}), (\ref{FFFF2}). The proof is completed.
\end{proof}

\end{document}